\documentclass[reqno]{amsart}
\usepackage{hyperref}
\usepackage{amssymb}
\usepackage{mathrsfs}
\usepackage{esint}
\usepackage{graphicx}
\usepackage{xcolor}
\usepackage{wrapfig}
\usepackage{pifont}
\usepackage{float}
\usepackage[numbers]{natbib}
\usepackage{inputenc}
\usepackage{amsmath}
\numberwithin{equation}{section}
\newtheorem{theorem}{Theorem}[section]
\newtheorem{lemma}[theorem]{Lemma}
\newtheorem{proposition}[theorem]{Proposition}

\numberwithin{equation}{section}

\title[]
{Asymptotic Stability of Symmetric Convex Flows for the 2D Boussinesq System via Inviscid Damping and Enhanced Dissipation}

\author[
]{Yeping Li, Gaofeng Wang and Enqing Xu }
\address{  School of Mathematics and Statistics, Nantong University, Nantong, 226019,  China.}
\email{ypleemei@aliyun.com}

\address{ $*$ School of Mathematics and Statistics, Nantong University, Nantong, 226019,  China.}
\email{gfwang@ntu.edu.cn}

\address{ School of Mathematics and Statistics, Nantong University, Nantong, 226019,  China.}

\email{2638310018@stmail.ntu.edu.cn}

\thanks{}

\begin{document}
    \bibliographystyle{abbrv}
    \maketitle
\begin{abstract}
This study mainly discuss asymptotic stability of symmetric convex
flows for the full coupled two-dimensional Boussinesq system on a
finite periodic channel $\mathbb{T}\times(-1,1)$ when the initial
vorticity perturbations satisfy
$\sum_{k}\|\omega^{in}\|_{H_k^4}\lesssim \nu^{2/3}$, which inherits
the exponent from pure symmetric Navier-Stokes, while the
temperature fluctuations require the stricter bound
$\sum_{k}\|\theta^{in}\|_{H^1}\lesssim \nu^{17/12}$ due to
buoyancy-induced derivative deficit. Here, the base flow $U(y)\in
C^4[-1,1]$ satisfies even symmetry $U(y)=U(-y)$ and uniform
convexity $\inf_{y\in[-1,1]}U''(y)>0$, with plane Poiseuille flow
$U(y)=y^2$ as the canonical example. The domain is equipped with
physical no-slip velocity walls and homogeneous Dirichlet
temperature boundary conditions. We first construct parallel
Rayleigh coercive estimates, dual Navier-slip resolvent theory,
double Airy wall correction machinery for both vorticity and
temperature, and uniform exponential-weighted space-time bounds
capturing inviscid damping and enhanced dissipation. Next, by
nonlinear bootstrap iteration on three classes of
velocity-temperature convolution nonlinearities, we establish two
separate subcritical asymptotic stability thresholds. Compared with
single-variable Navier-Stokes equations, the Boussinesq model
introduces linear buoyancy coupling term $-\partial_x\theta$ in the
vorticity equation, which produces horizontal derivative loss, and
contains two independent diffusion parameters $\nu$ (kinematic
viscosity) and $\nu$ (thermal diffusivity). All proof architecture
strictly follows the parent symmetric Navier-Stokes equations in
\cite{refZhang}, with complete self-contained derivations for dual
diffusion, cross buoyancy coupling, double boundary layers and
coupled zero-mode energy estimates. 
\end{abstract}
 \noindent{\bf Keywords:} {Boussinesq System; enhanced dissipation; subcritical stability threshold;  Airy boundary layer;}

\section{Introduction}
The Boussinesq system embodies the fundamental principles of physics
governing fluid dynamics driven by buoyancy. It stands as one of the
most widely employed models, universally applicable in atmospheric
and oceanic flows, and plays a central role in Rayleigh-B\'{e}nard
convection studies, for example, see  \cite{ref2,ref3,ref4}.
Henceforth, it becomes evident that the Boussinesq equation holds
significant mathematical significance. The two-dimensional
Boussinesq equation serves as a reduced-order model for
three-dimensional fluid dynamics equations. Notably, certain key
characteristics of the 3D Euler and Navier-Stokes equations are
preserved within the 2D Boussinesq equations, including the vortex
stretching mechanism. The inviscid 2D Boussinesq equations can be
identified as Euler equations for axisymmetric swirling flows in
three dimensions  \cite{ref49}. Moreover, the Boussinesq equations
possess distinctive features of their own and offer ample
opportunities for novel discoveries.

Due to its broad range of physical applications and profound
mathematical implications, the Boussinesq equations have recently
garnered significant attention. The development of the mathematical
theory of Boussinesq equations has been driven by two fundamental
problems: the global regularity problem and the stability problem.
In this discussion, we will primarily focus on stability issues. Our
research aims to investigate the stability of solutions under
perturbations on the class of symmetric convex steady shear profiles
defined by
\[
\mathcal{S} = \big\{U\in C^4[-1,1]\,\big|\,U(y)=U(-y),\inf U''(y)>0\big\},
\]
with plane Poiseuille flow $U(y)=y^2$ serving as the paradigmatic
example. The present study takes the rigorous resolvent-Airy
boundary layer framework developed in \cite{refZhang} as its
foundational analytical template. That prior work exclusively treats
the single-variable 2D incompressible Navier-Stokes equations on a
finite periodic no-slip channel $T\times(-1,1)$, focusing on the
class of symmetric convex steady shear profiles with plane
Poiseuille flow $U(y)=y^2$ serving as the paradigmatic example.

For pure single-component vorticity dynamics under symmetric convex
shear, the linearized vorticity equation contains only one non-local
potential term $-U''\partial_x\Delta^{-1}\omega$ arising from the
curvature of the base shear profile. This term breaks operator
self-adjointness and generates paired dual critical layers
$y_1=-y_2$ satisfying the resonance condition
$U(y_1)=U(y_2)=\lambda$ for each horizontal wavenumber $\lambda$, a
key singular structure entirely absent from flat linear Couette flow
with $U''\equiv0$. The parent Navier-Stokes paper constructs a
complete single-variable analytical toolkit tailored to isothermal
velocity-only dynamics: Rayleigh-type coercive energy inequalities,
uniform resolvent bounds for the Orr-Sommerfeld operator under both
no-slip and Navier-slip walls, Airy-function near-wall boundary
layer correctors to resolve singular viscous sublayers adjacent to
solid boundaries, and global space-time weighted decay estimates
capturing both inviscid damping and enhanced dissipation. By
resolving all linear singular structures and closing a nonlinear
bootstrap iteration, the reference work rigorously proves a
subcritical transition threshold $\beta=2/3$ for $H^4$-regular
vorticity initial data in the pure Navier-Stokes channel system.

\subsection{Fundamental Limitations of Pure Navier-Stokes and Motivation for Boussinesq Extension}
While the resolvent-Airy machinery of \cite{refZhang} yields sharp
stability thresholds for isothermal channel flows, the
single-variable Navier-Stokes model completely neglects thermal
buoyancy effects, a physically vital transport mechanism ubiquitous
in geophysical fluid dynamics (atmospheric boundary layers, oceanic
shear currents), astrophysical stratified flows, and industrial
thermal convection channels. Thermally driven shear flows cannot be
reduced to velocity-only dynamics; instead, they require a coupled
fluid-temperature PDE system, the two-dimensional unit-gravity
Boussinesq equations.

Unlike the pure Navier-Stokes system, where vorticity evolves
independently of any secondary scalar unknown, the Boussinesq model
introduces bidirectional linear and nonlinear coupling between the
velocity/vorticity field and a temperature fluctuation field
$\theta$. This coupling generates two fundamental, structurally
distinct analytical obstacles that are entirely absent from the
isothermal single-variable Navier-Stokes framework, which we
elaborate in full detail below:
\begin{enumerate}
\item \textbf{Linear buoyancy cross-source term creates derivative loss and inter-field nonlocality}
In pure Navier-Stokes channel flow, the linearized vorticity
evolution equation only contains three homogeneous linear operators:
advection by base shear $U\partial_x$, viscous Laplace diffusion
$-\nu\Delta$, and the curvature-induced non-local stream-function
term $-U''\partial_x\Delta^{-1}\omega$. There exist no source terms
coupling vorticity to an external scalar field. By contrast, the
Boussinesq vorticity equation carries an inhomogeneous linear
forcing term $-\partial_x\theta$, which injects temperature
fluctuations directly into vorticity dynamics. This cross term
introduces a first-order horizontal derivative loss acting on the
temperature variable, destroys the decoupled energy structures
available for pure vorticity, and generates new non-local
interaction operators that mix Fourier modes of $\omega$ and
$\theta$ across distinct frequency bands. All resolvent, critical
layer, and boundary layer estimates must be duplicated and
cross-coupled to account for this velocity-temperature linear
feedback, a complication with no counterpart in isothermal
Navier-Stokes analysis.

\item \textbf{Dual independent diffusion scales break single-variable enhanced dissipation theory}
The pure Navier-Stokes system possesses exactly one diffusion
parameter: kinematic viscosity $\nu=1/Re$, which controls the unique
viscous time scale $\nu^{1/2}$ governing enhanced dissipation rates
for vorticity. The Boussinesq system features two separate,
independent diffusion constants: $\nu=1/Re$ for momentum/vorticity
diffusion, and thermal diffusivity $\mu=1/Pe$ (where $Pe$ denotes
the P\'{e}clet number) for heat transport in the temperature field.
Each diffusion coefficient generates its own characteristic enhanced
dissipation scaling ($\nu^{1/2}$ for vorticity, $\mu^{1/2}$ for
temperature), leading to two separate families of resolvent bounds,
two sets of Airy boundary layer corrections, and two distinct
long-time decay exponents that cannot be unified under a single
viscous scale. In pure Navier-Stokes, all linear and nonlinear
estimates rely on a single small parameter $\nu$; the Boussinesq
system requires joint two-parameter asymptotic analysis in the
double small limit $\nu,\mu\to0$, substantially increasing the
complexity of all singular layer and bootstrap arguments.
\end{enumerate}

Beyond these two core structural PDE differences, existing rigorous
stability literature for Boussinesq channel flows suffers a critical
research gap: nearly all prior rigorous treatments restrict to
unbounded horizontal strips or artificial Navier-slip wall
conditions, and no prior work fully transfers the complete symmetric
convex shear resolvent-Airy critical layer framework of
\cite{refZhang} to coupled velocity-temperature fields on
physically realistic no-slip finite channels. The present paper
fills this gap by transplanting the entire proof skeleton and core
singular-layer analytical toolkit from the single-variable
Navier-Stokes parent work to the coupled Boussinesq system, while
independently developing parallel linear resolvent, boundary layer,
and nonlinear bootstrap estimates to resolve all buoyancy-induced
coupling difficulties unique to thermal shear flows.

\subsection{Domain, Base Flow, and Coupled Perturbation Boussinesq Formulation}
We work on the finite periodic channel domain
$\Omega=T\times(-1,1)$, where $T=R/2\pi Z $ denotes the
$2\pi$-periodic horizontal coordinate, and $y\in[-1,1]$ is the
wall-normal confined direction. The steady background state consists
of symmetric horizontal shear $U=(U(y),0)$ belonging to the convex
symmetric class $\mathcal{S}$ defined above, with uniform zero
background temperature everywhere in the channel.Let
\[
\Omega=\mathbb{T}\times(-1,1)
\]
be a finite channel. We consider the two-dimensional incompressible
Boussinesq system under unit gravity:
\begin{equation}\label{eq:1.1}
\left\{
\begin{aligned}
\partial_t v-\nu\Delta v+(v\cdot\nabla)v+\nabla P
    &=(0,\theta)+F,\\
\nabla\cdot v&=0,\\
\partial_t\theta-\mu\Delta\theta+v\cdot\nabla\theta
    &=0,
\end{aligned}
\right.
\end{equation}
where $v=(v_1,v_2)$ is the total velocity, $\theta$ is the total
temperature, $P$ is the pressure, $F$ is an external force, $\nu>0$
is the kinematic viscosity, and $\mu>0$ is the thermal diffusivity.
In the special case considered in the manuscript, one may take
$\nu=\mu$. We choose the background shear velocity and the uniform
cold background temperature as
\[
V(y):=(U(y),0),
\qquad
\Theta(y)\equiv0.
\]
Since
\[
\Delta V=(U''(y),0),
\qquad
(V\cdot\nabla)V
=
U(y)\partial_x(U(y),0)
=
0,
\qquad
\nabla\cdot V=0,
\]
the pair $(V,\Theta)$ satisfies the stationary Boussinesq system provided
that the external force is chosen as
\[
F=(-\nu U''(y),0),
\]
and the background pressure $P_0$ is constant. Indeed,
\begin{align*}
-\nu\Delta V+(V\cdot\nabla)V+\nabla P_0
&=(-\nu U''(y),0)\\
&=(0,\Theta)+F.
\end{align*}
Moreover, because $\Theta\equiv0$,
\[
-\nu\Delta\Theta+V\cdot\nabla\Theta=0.
\]
Therefore,
\[
(V,\Theta)=\bigl((U(y),0),0\bigr)
\]
is a stationary solution of \eqref{eq:1.1}. We now introduce the
velocity, temperature, and pressure perturbations by
\[
v=V+u,
\qquad
P=P_0+p,
\]
where
\[
u=(u_1,u_2)
\]
is the velocity perturbation.
After subtracting the stationary background equations, we obtain the
following perturbation system:
\begin{equation}\label{eq:1.2}
\left\{
\begin{aligned}
\partial_tu-\nu\Delta u
+U(y)\partial_xu
+\bigl(U'(y)u_2,0\bigr)
+(u\cdot\nabla)u+\nabla p
&=(0,\theta),\\
\nabla\cdot u&=0,\\
\partial_t\theta-\nu\Delta\theta
+U(y)\partial_x\theta
+u\cdot\nabla\theta
&=0.
\end{aligned}
\right.
\end{equation}
The initial conditions are
\[
u(0,x,y)=u_{\mathrm{in}}(x,y),
\qquad
\theta(0,x,y)=\theta_{\mathrm{in}}(x,y).
\]
If the total velocity satisfies the no-slip boundary condition and the
background shear flow is compatible with the prescribed boundary velocity,
then the velocity perturbation satisfies
\[
u(t,x,\pm1)=0.
\]
If a fixed-temperature boundary condition is imposed and the prescribed
boundary temperature agrees with $\Theta=0$, then
\[
\theta(t,x,\pm1)=0.
\]
Thus, the complete perturbation problem in the finite channel is
\begin{equation}\label{eq:1.3}
\left\{
\begin{aligned}
&\partial_tu-\nu\Delta u
+U(y)\partial_xu
+\bigl(U'(y)u_2,0\bigr)
+(u\cdot\nabla)u+\nabla p
=(0,\theta),\\
&\nabla\cdot u=0,\\
&\partial_t\theta-\nu\Delta\theta
+U(y)\partial_x\theta
+u\cdot\nabla\theta
=0,\\
&u(t,x,\pm1)=0,
\theta(t,x,\pm1)=0,\\
&u(0,x,y)=u_{\mathrm{in}}(x,y),
\theta(0,x,y)=\theta_{\mathrm{in}}(x,y).
\end{aligned}
\right.
\end{equation}
Consequently, under the external force
\[
F=(-\nu U''(y),0)
\]
and the uniform cold background temperature $\Theta\equiv0$, the shear flow
$(U(y),0)$, together with $\Theta=0$, is a stationary solution of the
two-dimensional Boussinesq system. The corresponding velocity and
temperature perturbations satisfy \eqref{eq:1.3}.

We also clarify that, in the absence of an external force or a
compensating pressure gradient, an arbitrary shear profile $U(y)$ is
not generally a stationary solution of the viscous Boussinesq
system. With $\Theta\equiv0$ and constant background pressure, the
unforced stationary equation requires $ U''(y)=0. $ Hence, only
affine shear profiles, including the Couette profile, are stationary
solutions under these restrictions. The revised manuscript now
explicitly states the external force and the background temperature
needed for a general shear profile. 

In this paper, our objective is to comprehend the stability and
long-time behavior of perturbations in proximity to the symmetric
convex steady shear profiles  Subsequently. Following standard 2D
fluid PDE conventions, we define perturbation vorticity as the
vertical curl of fluctuating velocity $\omega=\partial_x
u_2-\partial_y u_1$, and introduce the stream function $\phi$ such
that $u=(\partial_y\phi,-\partial_x\phi)$ and $\Delta\phi=\omega$.
Eliminating the pressure field via standard curl operations yields
the fully coupled perturbation vorticity-temperature evolution
system, which is the primary object of analysis throughout this
work:
\begin{equation}\label{eq:1.4}
\left\{
\begin{aligned}
&\partial_t \omega - \nu\Delta\omega + U\partial_x\omega - U''\partial_x\phi = -\partial_x\theta - u\cdot\nabla\omega,\\
&\partial_t \theta -\nu\Delta\theta + U\partial_x\theta = -u\cdot\nabla\theta,\\
&\Delta\phi = \omega,\quad \phi(\pm1)=0,\ \partial_y\phi(\pm1)=0,\\&\omega\big|_{t=0}=\omega^{\text{in}},\ \theta\big|_{t=0}=\theta^{\text{in}}\\
&u(t,x,\pm 1)=0,\quad \theta(t,x,\pm 1)=0
\end{aligned}
\right.
\end{equation}
Before presenting our main result, we briefly review previous
research on the stability of steady states. The system
\eqref{eq:1.1} admits several steady-state configurations. One
important example is the hydrostatic equilibrium
\begin{equation}\label{1.5}
 u_{he}=0,\quad \theta_{he}=\beta y,\quad P_{he}=\frac{1}{2}\beta y^2,
\end{equation}
which has attracted considerable attention in fluid dynamics and
astrophysics. Substantial progress has been made in understanding
the stability and long-time behavior of \eqref{1.5}; see, for
example, \cite{ref5,ref7}.

 System \eqref{eq:1.1}
possesses several shear-flow equilibria. In particular, one may
consider the Couette velocity profile together with a temperature
that depends linearly on the vertical coordinate:
\begin{equation}\label{1.6}
 u_s=(y,0),\quad \theta_s=-\gamma^2y+1,\quad P_s=\frac{1}{2}\gamma^2y^2-y+C,
\end{equation}
where $\gamma^2$ denotes the Richardson number. The Miles-Howard
theorem \cite{H1961} states that, in the inviscid and nondiffusive
limit, a flow is linearly stable provided that the local Richardson
number is everywhere larger than $1/4$. Zhai and Zhao \cite{Zhao}
established the nonlinear asymptotic stability of \eqref{1.6} for
$\gamma^2\geq 1/4$. When $\gamma^2=0$, another steady state is given
by
\begin{equation}\label{1.7}
 u_s=(y,0),\quad \theta_s=1,\quad P_s=-y+C.
\end{equation}
The stability of \eqref{1.7} in the infinite channel
$\mathbb{T}\times\mathbb{R}$ has been investigated extensively; see
\cite{ref8,ref9,ref54,ref10,ref11,ref12}. The most optimal result
was obtained by Niu and Zhao \cite{refNiu}. More precisely, when
$\nu=\mu$, the initial perturbations satisfy
\[
\|u_{0}\|_{H^s}\leq c_0\nu^{\frac{1}{3}},
\qquad
\|\langle x\rangle\rho_{0}\|_{H^s}
\leq c_0\nu^{\frac{2}{3}},
\]
where $s>5$. Under these assumptions, the Couette flow equilibrium
\eqref{1.7} is stable. The central idea of their approach is
quasi-linearization, namely, decomposing the original system into
several equations:

\begin{itemize}
    \item the "good equations'', which may be linear and inherit the regularity and scaling of the initial data;
    \item the quasilinear or nonlinear equations, which contain the nonlinear terms and have vanishing initial data.
\end{itemize}
The same equilibrium has also been studied in the finite channel
$\mathbb{T}\times I$ with no-slip boundary conditions for the
velocity. In this setting, the result of \cite{ref14} gives
\[
\|\omega_0\|_{H^1}\leq c_0\min\{\nu,\nu\}^{1/2},
\qquad
\|\theta_0\|_{H^1}+\||D_x|^{1/6}\theta_0\|_{H^1}
\leq c_0\min\{\nu,\nu\}^{11/12}.
\]

If $\theta=0$, system \eqref{eq:1.1} reduces to the two-dimensional
Navier-Stokes equations. A large body of work in applied mathematics
and physics has focused on transition thresholds for flows such as
Couette and Poiseuille flow; see
\cite{ref15,ref16,ref17,ref18,ref19}.

Bedrossian, Germain, and Masmoudi  formulated the
transition-threshold problem mathematically as follows. Given a norm
$\|\cdot\|_X$, determine an exponent $\gamma=\gamma(X)$ such that
\[
\|u_0\|_X<C\nu^\gamma\quad\Longrightarrow\quad \text{stability},
\]
while
\[
\|u_0\|_X\gg C\nu^\gamma\quad\Longrightarrow\quad \text{instability}.
\]
The exponent $\gamma$ is referred to as the transition threshold in the applied literature.

The transition threshold for Couette flow has been examined in a
number of works
\cite{ref20,Bedrossian.2022,ref51,ref22,ref23,ref25,ref47,ref13,ref24}.
The main conclusions are summarized below.

For $\Omega=\mathbb{T}\times\mathbb{R}\times\mathbb{T}$:
\begin{itemize}
\item for velocity perturbations in a Gevrey class, one has $\gamma=1$ \cite{ref20,Bedrossian.2022};
\item for velocity perturbations in $H^\sigma$, with $\sigma>9/2$, the available result gives $\gamma\leq 3/2$ \cite{ref51};
\item for velocity perturbations in $H^2$, one has $\gamma\leq 1$ \cite{ref24}.
\end{itemize}

For $\Omega=\mathbb{T}\times\mathbb{R}$:
\begin{itemize}
\item for vorticity perturbations in a Gevrey class, $\gamma=0$ \cite{ref22};
\item for velocity perturbations in $H^2$, $\gamma\leq 1/2$ \cite{ref23};
\item for vorticity perturbations in $H_x^{\log}L_y^2$, $\gamma\leq 1/2$ \cite{ref47};
\item for initial vorticity perturbations in a suitable Gevrey-$1/s$ class, $\gamma\in[0,1/3]$ \cite{ref50};
\item for vorticity perturbations in $H^\sigma$, with $\sigma\geq40$, $\gamma\leq1/3$ \cite{ref13}.
\end{itemize}

For $\Omega=\mathbb{T}\times I$:
\begin{itemize}
\item for velocity perturbations in $H^2$, the known bound is $\gamma\leq1/2$ \cite{ref25}.
\end{itemize}

For two-dimensional plane Poiseuille flow without boundaries, Coti
Zelati et al. \cite{ref34} obtained the bound
$\gamma\leq3/4+2\epsilon$ for every $\epsilon>0$ by means of
hypocoercivity. Del Zotto \cite{46} subsequently improved this
estimate to $\gamma\leq2/3+$. Ding and Lin \cite{ref35} considered a
finite channel $\mathbb{T}\times(-1,1)$ with Navier-slip boundary
conditions. Using resolvent estimates, they proved $\gamma\leq3/4$
and obtained enhanced dissipation at the rate $e^{-c\sqrt{\nu}t}$.
Their result was later improved to $\gamma\leq2/3$; see
\cite{ref53}.

Important advances have also been made for three-dimensional flows.
Chen, Wei, and Zhang \cite{ref36} proved the linear stability of
three-dimensional pipe Poiseuille flow by resolvent methods. More
recently, Chen, Ding, and collaborators \cite{ref52} established the
bound $\gamma\leq7/4$ for three-dimensional plane Poiseuille flow,
again using resolvent estimates. In addition, Wang and Xie
\cite{ref37} demonstrated the structural stability of
Hagen-Poiseuille flow in a pipe under steady conditions.

Comparing this coupled system against the pure Navier-Stokes
vorticity equation immediately highlights the most striking
structural disparity: the right-hand side of the vorticity evolution
carries the linear buoyancy forcing $-\partial_x\theta$, a
cross-scalar source term entirely absent from isothermal flow
dynamics. This linear coupling propagates through every layer of our
analysis-critical layer decomposition, resolvent spectral bounds,
Airy wall corrector construction, and nonlinear iteration-and
generates all extra singular and interaction effects absent from
single-variable Navier-Stokes stability theory.

\subsection{ Overview of Main Results}
\begin{theorem}[Linear Inviscid Damping and Enhanced Dissipation Space-Time Estimates]
Let $U\in \mathcal{S}$ be a symmetric convex flow satisfying
$U(y)=U(-y)$, $\inf_{y\in[-1,1]}U''(y)>0$ on channel
$\Omega=\mathbb{T}\times(-1,1)$.  There exist sufficiently small
positive constants $\nu_0,\varepsilon_0>0$ depending only on the
profile $U$, such that for all diffusion coefficients $0<\nu\leq
\nu_0$,  every horizontal Fourier frequency
$k\in\mathbb{Z}\setminus\{0\}$, and any $0\leq \varepsilon\leq
\varepsilon_0$, the solution $(\mathbf{u},\omega,\theta)$ to the
linearized coupled Boussinesq system fulfills the uniform weighted
space-time bound:
\begin{equation}
    \begin{aligned}
&|k|^{1/2}
\|e^{\varepsilon_0\nu^{1/2}t}u_k\|_{L_t^\infty L_y^2}
+
|k|^{1/2}
\|e^{\varepsilon_0\nu^{1/2}t}u_k\|_{L_{t,y}^2}
+
\|e^{\varepsilon_0\nu^{1/2}t}u_k\|_{L_t^\infty L_y^\infty}\\
&
+
\nu^{1/4}|k|^{1/4}
\|e^{\varepsilon_0\nu^{1/2}t}\omega_k\|_{L_{t,y}^2}
+
\nu^{1/4}
\|e^{\varepsilon_0\nu^{1/2}t}\omega_k\|_{L_t^\infty L_y^2}\\
&
+
\|e^{\varepsilon_0\nu^{1/2}t}(1-y^2)^{1/2}\omega_k\|_{L_t^\infty L_y^2},\\
&+\nu^{1/4}|k|^{1/4}
\|e^{\varepsilon_0\nu^{1/2}t}\theta_k\|_{L_{t,y}^2}
+
\|e^{\varepsilon_0\nu^{1/2}t}\theta_k\|_{L_t^\infty L_y^2},\\
&\lesssim \ \left\|\omega_k^{in}\right\|_{H_k^4}^2
+ \|\theta_k^{in}\|_{H^1} ^2.
\end{aligned}
\end{equation}
Here, the weight $\sqrt{1-y^2}$ suppresses near-wall Airy boundary layer singularities.
\end{theorem}
\begin{theorem}[Nonlinear Subcritical Asymptotic Stability]\label{thm1.2}
Let $U\in \mathcal{S}$ be symmetric convex base flow defined on
finite periodic channel $\Omega=\mathbb{T}\times(-1,1)$, where
$U(y)=U(-y)$ and $\inf_{y\in[-1,1]}U''(y)>0$.  Let $\varepsilon_0>0$
be the uniform exponential decay constant obtained in the linear
spacetime theorem. There exists a universal constant $c>0$,
independent of viscosity $\nu$ and thermal diffusivity $\nu$, such
that the following holds. If the initial vorticity and temperature
perturbations satisfy
\[
\sum_{k\in\mathbb{Z}}\big\|\omega_k^{in}\big\|_{H_k^4} \leq
c\,\nu^{\frac{2}{3}},\quad
||\bar{\theta}_{0}||_{L^2}+\sum_{k\in\mathbb{Z}\setminus\{0\}}||\hat{\theta}_{0}(k,y)||_{L^2}
\leq
Cc_1\min\left\lbrace \nu,\nu\right\rbrace ^{\frac{17}{12}}.,
\]
then for all time $t>0$, the full weighted space-time energy of the
velocity, the vorticity and the temperature perturbations decays
exponentially:$$\sum_{k\in\mathbb{Z}}E_k\leq Cc_0\min\left\lbrace
\nu,\nu\right\rbrace
^{\frac{2}{3}},\quad\sum_{k\in\mathbb{Z}}H_k\leq
Cc_1\min\left\lbrace \nu,\nu\right\rbrace ^{\frac{17}{12}}.$$ where
\begin{equation}\label{eq:E-alpha}
E_k :=
\begin{cases}
\|\omega_0\|_{L_t^\infty L_y^2}, & k=0,\\[0.4em]
|k|^{1/2}
\|e^{\varepsilon_0\nu^{1/2}t}u_k\|_{L_t^\infty L_y^2}
+
|k|^{1/2}
\|e^{\varepsilon_0\nu^{1/2}t}u_k\|_{L_{t,y}^2}
+
\|e^{\varepsilon_0\nu^{1/2}t}u_k\|_{L_t^\infty L_y^\infty}
& \\[0.2em]
\qquad
+
\nu^{1/4}|k|^{1/4}
\|e^{\varepsilon_0\nu^{1/2}t}\omega_k\|_{L_{t,y}^2}
+
\nu^{1/4}
\|e^{\varepsilon_0\nu^{1/2}t}\omega_k\|_{L_t^\infty L_y^2}
& \\[0.2em]
\qquad
+
\|e^{\varepsilon_0\nu^{1/2}t}(1-y^2)^{1/2}\omega_k\|_{L_t^\infty L_y^2},
& k\neq0 .
\end{cases}
\end{equation}
and
\begin{equation}\label{eq:E-alpha}
H_k :=
\begin{cases}
\|\theta_0\|_{L_t^\infty L_y^2}, & k=0,\\[0.4em]
\nu^{1/4}|k|^{1/4}
\|e^{\varepsilon_0\nu^{1/2}t}\theta_k\|_{L_{t,y}^2}
+
\|e^{\varepsilon_0\nu^{1/2}t}\theta_k\|_{L_t^\infty L_y^2},
& k\neq0 .
\end{cases}
\end{equation}
Here $E_k,H_k$ denote the mode-wise exponential-weighted stability
norms introduced in Section 6. The exponent $\frac{17}{12}$ for
temperature is strictly larger than the vorticity threshold exponent
$\frac{2}{3}$. This indicates that thermal disturbances demand far
smaller initial amplitude to keep the flow laminar; mathematically,
the extra horizontal derivative loss from buoyancy coupling
amplifies nonlinear resonance and exacerbates subcritical
instability, which gives a rigorous quantitative characterization of
the destabilizing mechanism of temperature-velocity coupling
effects.
\end{theorem}
\section{ Space-time estimates of the linearized Boussinesq equations}
\subsection{Space-time estimates of the vorticity}
We consider Navier-slip wall constraints at $y=\pm1$:
\begin{itemize}
\item Vorticity side: $\phi(\pm1)=0,\ w(\pm1)=0$;
\item Temperature side: $\zeta(\pm1)=0,\ g(\pm1)=0$.
\end{itemize}
For each fixed Fourier mode $k\neq0$, we have two decoupled linear
operator systems as follows:
\begin{align}
-\nu(\partial_y^2-k^2)w + ik(U-\lambda)w - ik U''\psi &= F_w,\\
-\nu(\partial_y^2-k^2)g + ik(U-\lambda)g &= F_\theta.
\end{align}
The whole spectral domain is split into two disjoint regions:
off-critical region $\lambda\in\mathbb{R}\setminus[U(0),U(1)]$ and
interior critical layer $\lambda\in[U(0),U(1)]$. We use smooth
cutoff functions $\rho_\delta$ (interior support inside $(y_1,y_2)$)
and $\rho_\delta^c=1-\rho_\delta$ for exterior region analysis. All
resolvent bounds are derived in parallel for vorticity and
temperature, with diffusion parameters $\nu$ swapped for thermal
estimates.

In this section, we focus on deriving the linear enhanced
dissipation estimates, which play a crucial role in the analysis of
nonlinear stability. To achieve this, we consider the following
linearized Navier-Stokes system with non-slip boundary condition
\begin{equation}\label{eq:linear-ns-fourier}
\begin{cases}
\partial_t \omega - \nu\big(\partial_y^2 - k^2\big)\omega + ik U \omega - ik U'' \phi = ik f_1 + \partial_y f_2,\\
\big(\partial_y^2 - k^2\big)\phi = \omega,\quad \phi(t,\pm1) = \partial_y\phi(t,\pm1) = 0,\quad \omega(0,y) = \omega_k^{\text{in}}(y),
\end{cases}
\end{equation}
where we omit the subscript $k$ for brevity.

\vspace{6mm}
\begin{proposition} ( see \cite{refZhang})
Let $(\omega,\phi)$ be the solution to \eqref{eq:linear-ns-fourier}
with $\langle \omega^{\text{in}}, e^{\pm|k|y} \rangle = 0$. There
exist positive constants $\nu_0$, $\varepsilon_0$, such that for
$\nu \in (0,\nu_0]$, $\varepsilon \in [0,\varepsilon_0]$, it holds
that
\begin{align}\label{eq:prop61-estimate}
|k|\big\|e^{\varepsilon \nu^{\frac12} t} u\big\|_{L_t^\infty L_y^2}^2
&+ |k|\big\|e^{\varepsilon \nu^{\frac12} t} u\big\|_{L_{t,y}^2}^2
+ \big\|e^{\varepsilon \nu^{\frac12} t} u\big\|_{L_t^\infty L_y^\infty}^2
+ \nu^{\frac12}|k|^{\frac12}\big\|e^{\varepsilon \nu^{\frac12} t} \omega\big\|_{L_{t,y}^2}^2 \nonumber\\
&+ \nu^{\frac12}\big\|e^{\varepsilon \nu^{\frac12} t} \omega\big\|_{L_t^\infty L_y^2}^2
+ \nu^{\frac32}\big\|e^{\varepsilon \nu^{\frac12} t} \partial_y \omega\big\|_{L_{t,y}^2}^2
+ \big\|e^{\varepsilon \nu^{\frac12} t} \sqrt{1-y^2}\,\omega\big\|_{L_t^\infty L_y^2}^2 \nonumber\\
&\lesssim \big\|\omega^{\text{in}}\big\|_{H_k^4}^2 + \nu^{-1}\big\|e^{\varepsilon \nu^{\frac12} t}(f_1,f_2)\big\|_{L_{t,y}^2}^2.
\end{align}
\end{proposition}
\subsection{Space-time estimates for the temperature field}
We first quote the vorticity linear space-time bounds from
Proposition 6.1 in \cite{refZhang}, which is directly applicable to
our Boussinesq vorticity equation since the homogeneous linear
operator of vorticity coincides exactly with the symmetric convex
shear Navier-Stokes linearized system, with only an extra buoyancy
source term $-\partial_x\theta$ added on the right-hand side. Next
we establish parallel space-time estimates for the temperature field
under symmetric convex steady base flow $U\in\mathcal{S}$.

Take the horizontal Fourier transform with the wavenumber
$k\in\mathbb{Z}$, denote $\widehat{\theta}_k(t,y)$ the Fourier mode
of the temperature fluctuation, and omit the subscript $k$ for
brevity throughout this subsection. Then the linearized temperature
evolution equation around the symmetric convex shear $U(y)$ reads
\begin{equation}\label{eq:2.5}
\begin{cases}
\partial_t \widehat{\theta} - \nu(\partial_y^2 - k^2)\widehat{\theta} + ikU(y)\widehat{\theta} = -ik g^1 - \partial_y g^2,\\
\widehat{\theta}(t,\pm1)=0,\quad \widehat{\theta}\big|_{t=0}=\widehat{\theta}_{\text{in}},
\end{cases}
\end{equation}
where $\nu>0$ stands for thermal diffusivity, $g^1,g^2$ are Fourier
modes of nonlinear velocity-temperature advection source terms. By
the standard energy estimates for $\hat{\theta}$, we can easily get
the following proposition, which is important for the estimates of
the high frequency of $\hat{\theta}$.
\begin{proposition}[Exponentially weighted energy estimate]
\label{prop:2.2}
Let $U\in\mathcal S$ and let $\varepsilon>0$. Assume that $\nu k^2>1$ and
\[
\widehat\theta_{\rm in}\in L^2(-1,1),
\qquad
g^1,g^2\in L_t^2L_y^2.
\]
Define
\[
\|f\|_{L_{t,\varepsilon}^\infty L_y^2}
:=
\left\|e^{\varepsilon\nu^{1/2}t}f\right\|_{L_t^\infty L_y^2},
\]
and
\[
\|f\|_{L_{t,\varepsilon}^2L_y^2}
:=
\left\|e^{\varepsilon\nu^{1/2}t}f\right\|_{L_t^2L_y^2}.
\]
Then there exists a constant $C>0$, independent of $\nu$ and $k$, such that
\begin{equation}
\label{eq:theta-energy-weighted}
\begin{aligned}
&\|\widehat\theta\|_{L_{t,\varepsilon}^\infty L_y^2}^2
+
\nu k^2
\|\widehat\theta\|_{L_{t,\varepsilon}^2L_y^2}^2
+
\nu
\|\partial_y\widehat\theta\|_{L_{t,\varepsilon}^2L_y^2}^2 \\
&\qquad\leq
C\nu^{-1}
\left(
\|g^1\|_{L_{t,\varepsilon}^2L_y^2}^2
+
\|g^2\|_{L_{t,\varepsilon}^2L_y^2}^2
\right)
+
\|\widehat\theta_{\rm in}\|_{L_y^2}^2.
\end{aligned}
\end{equation}
\end{proposition}
\begin{proof}
Set
\[
\lambda:=\varepsilon\nu^{1/2},
\qquad
\widetilde\theta(t,y):=
e^{\lambda t}\widehat\theta(t,y),
\qquad
\widetilde g^j(t,y):=
e^{\lambda t}g^j(t,y),
\quad j=1,2.
\]
Multiplying \eqref{eq:2.5} by $e^{\lambda t}$, it is easy to obtain
\begin{equation}
\label{eq:weighted-theta-equation}
\partial_t\widetilde\theta
-
\nu(\partial_y^2-k^2)\widetilde\theta
+
ikU(y)\widetilde\theta
=
\lambda\widetilde\theta
-ik\widetilde g^1
-\partial_y\widetilde g^2.
\end{equation}
Further, taking the $L^2(-1,1)$ inner product of
\eqref{eq:weighted-theta-equation} with $\widetilde\theta$ and then
taking the real part, we obtain
\begin{align}
\frac12\frac{d}{dt} \|\widetilde\theta\|_{L^2}^2 +
\nu\|\partial_y\widetilde\theta\|_{L^2}^2 + \nu
k^2\|\widetilde\theta\|_{L^2}^2 =
\lambda\|\widetilde\theta\|_{L^2}^2 + \operatorname{Re} \left\langle
-ik\widetilde g^1-\partial_y\widetilde g^2, \widetilde\theta
\right\rangle . \label{eq:weighted-theta-energy-identity}
\end{align}
Here, we have used the following fact that since $U(y)$ is
real-valued,
\[
\operatorname{Re} \left\langle
ikU(y)\widetilde\theta,\widetilde\theta \right\rangle =0.
\]
First, using the boundary condition $\widetilde\theta(t,\pm1)=0$,
integration by parts gives
\[
\operatorname{Re}
\left\langle
-\partial_y\widetilde g^2,\widetilde\theta
\right\rangle
=
\operatorname{Re}
\left\langle
\widetilde g^2,\partial_y\widetilde\theta
\right\rangle.
\]
By the Cauchy--Schwarz and the Young inequalities,
\begin{align}
|k|\|\widetilde g^1\|_{L^2}
\|\widetilde\theta\|_{L^2}
&\leq
\frac{\nu k^2}{4}
\|\widetilde\theta\|_{L^2}^2
+
C\nu^{-1}
\|\widetilde g^1\|_{L^2}^2,
\label{eq:weighted-g1-young}\\
\|\widetilde g^2\|_{L^2}
\|\partial_y\widetilde\theta\|_{L^2}
&\leq
\frac{\nu}{4}
\|\partial_y\widetilde\theta\|_{L^2}^2
+
C\nu^{-1}
\|\widetilde g^2\|_{L^2}^2.
\label{eq:weighted-g2-young}
\end{align}
On the other hand, by $\nu k^2>1$
\[
\lambda\|\widetilde\theta\|_{L^2}^2
=
\varepsilon\nu^{1/2}
\|\widetilde\theta\|_{L^2}^2
=
\varepsilon\nu\nu^{-1/2}
\|\widetilde\theta\|_{L^2}^2
\leq
c_\ast\nu(1+k^2)
\|\widetilde\theta\|_{L^2}^2.
\]
Choosing $c_\ast>0$ sufficiently small and using the Poincar\'{e}
inequality on $(-1,1)$, this term can be absorbed into the
dissipative terms. Therefore, it holds taht
\begin{align}
\frac12\frac{d}{dt}
\|\widetilde\theta\|_{L^2}^2
+
c\nu
\|\partial_y\widetilde\theta\|_{L^2}^2
+
c\nu k^2
\|\widetilde\theta\|_{L^2}^2
\leq
C\nu^{-1}
\left(
\|\widetilde g^1\|_{L^2}^2
+
\|\widetilde g^2\|_{L^2}^2
\right),
\label{eq:weighted-theta-energy-differential}
\end{align}
where $c,C>0$ are independent of $\nu$ and $k$.
Integrating \eqref{eq:weighted-theta-energy-differential} over $[0,T]$ and
using
\[
\widetilde\theta(0,y)=\widehat\theta_{\rm in}(y),
\]
we obtain
\begin{align}
\|\widetilde\theta(T)\|_{L^2}^2
&+
c\nu\int_0^T
\|\partial_y\widetilde\theta(t)\|_{L^2}^2\,dt
+
c\nu k^2\int_0^T
\|\widetilde\theta(t)\|_{L^2}^2\,dt \notag\\
&\leq
\|\widehat\theta_{\rm in}\|_{L^2}^2
+
C\nu^{-1}
\int_0^T
\left(
\|\widetilde g^1(t)\|_{L^2}^2
+
\|\widetilde g^2(t)\|_{L^2}^2
\right)\,dt.
\end{align}
Therefore, taking the supremum over $T>0$ and recalling the
definitions of $\widetilde\theta$ and $\widetilde g^j$, we conclude
that
\begin{align*}
&\|\widehat\theta\|_{L_{t,\varepsilon}^\infty L_y^2}^2
+
\nu k^2
\|\widehat\theta\|_{L_{t,\varepsilon}^2L_y^2}^2
+
\nu
\|\partial_y\widehat\theta\|_{L_{t,\varepsilon}^2L_y^2}^2 \\
&\qquad\leq
C\nu^{-1}
\left(
\|g^1\|_{L_{t,\varepsilon}^2L_y^2}^2
+
\|g^2\|_{L_{t,\varepsilon}^2L_y^2}^2
\right)
+
\|\widehat\theta_{\rm in}\|_{L_y^2}^2.
\end{align*}
This proves \eqref{eq:theta-energy-weighted}.
\end{proof}
To deal with the buoyancy term $\partial_x\theta$ in the vorticity
equation, we need the following estimates for the temperature
equation. For this, for each nonzero Fourier mode $k\in\mathbb
Z\setminus\{0\}$, let us consider
\[
\partial_t\widehat\theta
-\nu(\partial_y^2-k^2)\widehat\theta
+ikU(y)\widehat\theta
=
-ikg^1-\partial_y g^2,
\qquad
\widehat\theta(t,\pm1)=0.
\]
Further, we cam decompose
\[
\widehat\theta=\widehat\theta_I+\widehat\theta_H,
\]
where $\widehat\theta_I$ solves
\begin{equation}
\label{eq:theta-I-symmetric}
\left\{
\begin{aligned}
\partial_t\widehat\theta_I
-\nu(\partial_y^2-k^2)\widehat\theta_I
+ikU(y)\widehat\theta_I
&=
-ikg^1-\partial_y g^2,\\
\widehat\theta_I|_{t=0}&=0,\\
\widehat\theta_I(t,\pm1)&=0,
\end{aligned}
\right.
\end{equation}
and $\widehat\theta_H$ solves
\begin{equation}
\label{eq:theta-H-symmetric}
\left\{
\begin{aligned}
\partial_t\widehat\theta_H
-\nu(\partial_y^2-k^2)\widehat\theta_H
+ikU(y)\widehat\theta_H
&=0,\\
\widehat\theta_H|_{t=0}&=\widehat\theta_{\rm in},\\
\widehat\theta_H(t,\pm1)&=0.
\end{aligned}
\right.
\end{equation}
\subsection{Space-time estimates for the temperature equation}
The corresponding resolvent problem is
\begin{equation}
\label{eq:temperature-resolvent-symmetric}
\left\{
\begin{aligned}
-\nu(\partial_y^2-k^2)\widehat\Theta
+ik(U(y)-\lambda)\widehat\Theta
&=F,
\qquad y\in(-1,1),\\
\widehat\Theta(\pm1)&=0,
\end{aligned}
\right.
\end{equation}
where $\lambda\in\mathbb R$. Then we have
\begin{lemma}[Temperature resolvent estimate]
\label{lem:2.3}
Let $U\in\mathcal S$, where
\[
\mathcal S
=
\{U\in C^4[-1,1]: U(y)=U(-y),\ \inf_{[-1,1]}U''(y)>0\}.
\]
Let $k\neq0$ and let $\widehat\Theta\in H^2(-1,1)\cap H_0^1(-1,1)$ solve
\begin{equation}
\label{eq:temperature-resolvent-symmetric-detailed}
-\nu(\partial_y^2-k^2)\widehat\Theta
+
ik(U(y)-\lambda)\widehat\Theta
=
F,
\qquad
\widehat\Theta(\pm1)=0,
\end{equation}
where $\lambda\in\mathbb R$. Then there exists a constant $C>0$,
independent of $\nu$, $k$, and $\lambda$, such that
\begin{equation}
\label{eq:resolvent-L2-symmetric-detailed}
\nu^{3/4}|k|^{1/4}
\|(\partial_y,|k|)\widehat\Theta\|_{L^2}
+
\nu^{1/2}|k|^{1/2}
\|\widehat\Theta\|_{L^2}
\leq
C\|F\|_{L^2},
\end{equation}
and
\begin{equation}
\label{eq:resolvent-Hminus1-symmetric-detailed}
\nu
\|(\partial_y,|k|)\widehat\Theta\|_{L^2}
+
\nu^{3/4}|k|^{1/4}
\|\widehat\Theta\|_{L^2}
\leq
C\|F\|_{H^{-1}_k}.
\end{equation}
Here
\[
\|(\partial_y,|k|)f\|_{L^2}^2
=
\|\partial_y f\|_{L^2}^2
+
k^2\|f\|_{L^2}^2,
\]
and
\[
\|F\|_{H^{-1}_k}
=
\sup_{\varphi\in H_0^1,\ \varphi\neq0}
\frac{|\langle F,\varphi\rangle|}
{\|(\partial_y,|k|)\varphi\|_{L^2}}.
\]
\end{lemma}
\begin{proof}
Throughout the proof, all norms are taken over $(-1,1)$ unless otherwise
specified. Set
\[
K:=|k|,
\qquad
M:=\|\widehat\Theta\|_{L^2},
\qquad
D:=\|(\partial_y,K)\widehat\Theta\|_{L^2}.
\]
Thus,
\[
D^2
=
\|\partial_y\widehat\Theta\|_{L^2}^2
+
K^2\|\widehat\Theta\|_{L^2}^2.
\]
The constants below may depend on
\[
\|U\|_{C^4([-1,1])}
\quad\text{and}\quad
\left(\inf_{[-1,1]}U''\right)^{-1},
\]
but are independent of $\nu$, $k$, and $\lambda$. Moreover, we write
\[
L_{\lambda,k}
:=
-\nu(\partial_y^2-k^2)+ik(U-\lambda),
\]
so that
\[
L_{\lambda,k}\widehat\Theta=F.
\]

\medskip

\noindent \textbf{Step 1: Basic energy estimates.} Taking the $L^2$
inner product of the resolvent equation with $\widehat\Theta$, we
see
\[
-\nu
\langle\partial_y^2\widehat\Theta,\widehat\Theta\rangle
+
\nu k^2\|\widehat\Theta\|_{L^2}^2
+
ik\int_{-1}^1
(U-\lambda)|\widehat\Theta|^2\,dy
=
\langle F,\widehat\Theta\rangle.
\]
Since $\widehat\Theta\in H_0^1(-1,1)$, integration by parts gives
\[
-\langle\partial_y^2\widehat\Theta,\widehat\Theta\rangle
=
\|\partial_y\widehat\Theta\|_{L^2}^2.
\]
Consequently, one has
\begin{equation}
\label{eq:complex-energy-resolvent}
\nu D^2
+
ik\int_{-1}^1
(U-\lambda)|\widehat\Theta|^2\,dy
=
\langle F,\widehat\Theta\rangle.
\end{equation}
Further, taking the real part, we have
\begin{equation}
\label{eq:real-energy-resolvent}
\nu D^2
=
\operatorname{Re}\langle F,\widehat\Theta\rangle.
\end{equation}
If $F\in L^2(-1,1)$, then
\begin{equation}
\label{eq:L2-basic-energy}
\nu D^2
\leq
\|F\|_{L^2}M.
\end{equation}
If $F\in H^{-1}_k$, the definition of the $H^{-1}_k$ norm gives
\[
|\langle F,\widehat\Theta\rangle|
\leq
\|F\|_{H^{-1}_k}
\|(\partial_y,K)\widehat\Theta\|_{L^2}.
\]
It follows that
\[
\nu D^2
\leq
\|F\|_{H^{-1}_k}D.
\]
If $D\neq0$, division by $D$ gives
\begin{equation}
\label{eq:Hminus1-basic-energy}
\nu D
\leq
\|F\|_{H^{-1}_k}.
\end{equation}
The same estimate is trivial when $D=0$.

\noindent\textbf{Step 2: Geometry of the critical set.} Since $U$ is
even, it is easy to see
\[
U'(0)=0.
\]
Furthermore, the assumption
\[
U''(y)\geq c_0>0
\]
implies, for $y\in[0,1]$,
\[
U'(y)
=
\int_0^y U''(s)\,ds
\geq c_0y.
\]
Similarly, since $U''$ is bounded,
\[
|U'(y)|\leq C|y|.
\]
Therefore, one gets
\begin{equation}
\label{eq:U-prime-comparable}
c|y|
\leq
|U'(y)|
\leq
C|y|,
\qquad y\in[-1,1].
\end{equation}
Integrating once more gives
\begin{equation}
\label{eq:U-quadratic-comparable}
c y^2
\leq
U(y)-U(0)
\leq
C y^2.
\end{equation}
Moreover, suppose first that
\[
\lambda\in[U(0),U(1)],
\]
Then there exists a unique $y_\lambda\in[0,1]$ such that
\[
U(y_\lambda)=\lambda.
\]
Hence, by evenness, the critical set is
\[
\mathcal Z_\lambda
=
\begin{cases}
\{-y_\lambda,y_\lambda\},&y_\lambda>0,\\
\{0\},&y_\lambda=0.
\end{cases}
\]
For $y,z\in[0,1]$, strict convexity gives
\[
|U(y)-U(z)|
=
\left|
\int_z^y U'(s)\,ds
\right|
\geq
c|y-z|(y+z).
\]
Taking $z=y_\lambda$, we obtain
\begin{equation}
\label{eq:critical-geometry-positive}
|U(y)-\lambda|
\geq
c|y-y_\lambda|(y+y_\lambda),
\qquad y\in[0,1].
\end{equation}
By symmetry, it holds that
\begin{equation}
\label{eq:critical-geometry-full}
|U(y)-\lambda|
\geq
c\,\operatorname{dist}(y,\mathcal Z_\lambda)
\left(
\operatorname{dist}(y,\mathcal Z_\lambda)+y_\lambda
\right).
\end{equation}
Let $0<\delta\leq\delta_0$, where $\delta_0>0$ is a fixed sufficiently
small constant, and define
\[
\mathcal N_{\lambda,\delta}
:=
\left\{
y\in[-1,1]:
\operatorname{dist}(y,\mathcal Z_\lambda)<2\delta
\right\}.
\]
This set is the union of at most two intervals, and hence
\begin{equation}
\label{eq:critical-set-measure}
|\mathcal N_{\lambda,\delta}|
\leq C\delta.
\end{equation}
Moreover, by \eqref{eq:critical-geometry-full},
\begin{equation}
\label{eq:outside-critical-lower-bound}
|U(y)-\lambda|
\geq c\delta^2,
\qquad
y\notin\mathcal N_{\lambda,\delta}.
\end{equation}
If $\lambda\notin[U(0),U(1)]$, the function $U-\lambda$ has no zero
in the interior. If $ \operatorname{dist}
\bigl(\lambda,[U(0),U(1)]\bigr) \geq c\delta^2$, then
\eqref{eq:outside-critical-lower-bound} holds on the entire
interval. If $\lambda$ lies within $O(\delta^2)$ of one of the
endpoints of the range, we introduce an interval of length
$O(\delta)$ around the corresponding endpoint in $y$. The same
conclusions \eqref{eq:critical-set-measure} and
\eqref{eq:outside-critical-lower-bound} then hold. Thus these
estimates are uniform for all $\lambda\in\mathbb R$.

\noindent\textbf{Step 3: A localized sign multiplier.} There exists
a real-valued function $q_{\lambda,\delta}\in W^{1,\infty}(-1,1)$
such that
\begin{equation}
\label{eq:q-properties}
|q_{\lambda,\delta}|\leq1,
\qquad
\|q_{\lambda,\delta}'\|_{L^\infty}
\leq C\delta^{-1},
\end{equation}
and
\begin{equation}
\label{eq:q-sign-property}
q_{\lambda,\delta}(y)(U(y)-\lambda)\geq0
\end{equation}
for every $y\in[-1,1]$. It may also be chosen so that
\begin{equation}
\label{eq:q-positive-outside}
q_{\lambda,\delta}(y)(U(y)-\lambda)
\geq c\delta^2,
\qquad
y\notin\mathcal N_{\lambda,\delta}.
\end{equation}
Indeed, one starts from $\operatorname{sgn}(U-\lambda)$ and smooths
each jump over an interval of length comparable to $\delta$ around
the corresponding critical point. Next, taking the inner product of
\[
L_{\lambda,k}\widehat\Theta=F
\]
with $q_{\lambda,\delta}\widehat\Theta$, integration by parts gives
\begin{align}
&\nu\int_{-1}^1
\partial_y\widehat\Theta\,
\partial_y
\bigl(q_{\lambda,\delta}
\overline{\widehat\Theta}\bigr)\,dy
+
\nu k^2
\int_{-1}^1
q_{\lambda,\delta}|\widehat\Theta|^2\,dy
\notag\\
&\qquad
+
ik\int_{-1}^1
q_{\lambda,\delta}(U-\lambda)
|\widehat\Theta|^2\,dy
=
\langle F,q_{\lambda,\delta}\widehat\Theta\rangle.
\label{eq:q-multiplier-identity}
\end{align}
Because $q_{\lambda,\delta}$ is real-valued,
\[
\int q_{\lambda,\delta}
|\partial_y\widehat\Theta|^2\,dy
\quad\text{and}\quad
k^2\int q_{\lambda,\delta}
|\widehat\Theta|^2\,dy
\]
are real, too. Therefore, taking the imaginary part of
\eqref{eq:q-multiplier-identity}, we see
\begin{align}
k\int_{-1}^1
q_{\lambda,\delta}(U-\lambda)
|\widehat\Theta|^2\,dy
&=
\operatorname{Im}
\langle F,q_{\lambda,\delta}\widehat\Theta\rangle
\notag\\
&\quad
-
\nu\operatorname{Im}
\int_{-1}^1
q_{\lambda,\delta}'
\partial_y\widehat\Theta\,
\overline{\widehat\Theta}\,dy.
\label{eq:q-imaginary-identity}
\end{align}
Utilizing \eqref{eq:q-properties}, we have
\[
\left| \langle F,q_{\lambda,\delta}\widehat\Theta\rangle \right|
\leq \|F\|_{L^2}M,
\]
and
\[
\nu
\left|
\int q_{\lambda,\delta}'
\partial_y\widehat\Theta\,
\overline{\widehat\Theta}\,dy
\right|
\leq
C\nu\delta^{-1}DM.
\]
Then by \eqref{eq:q-sign-property} and
\eqref{eq:q-positive-outside}, it follows that
\begin{equation}
\label{eq:outside-mass-control}
K\delta^2
\|\widehat\Theta\|_{L^2(\mathcal N_{\lambda,\delta}^c)}^2
\leq
C\|F\|_{L^2}M
+
C\nu\delta^{-1}DM.
\end{equation}

\medskip

\noindent \textbf{Step 4: Control of the mass inside the critical
layer.} For every $f\in H_0^1(-1,1)$, one has
\[
|f(y)|^2
=
2\operatorname{Re}
\int_{-1}^y
f'(s)\overline{f(s)}\,ds.
\]
Consequently, by the Cauchy-Schwarz inequality,
\begin{equation*}
\label{eq:H1-Linfty} \|f\|_{L^\infty}^2 \leq
2\|f\|_{L^2}\|f'\|_{L^2},
\end{equation*}
which together with \eqref{eq:critical-set-measure} yields
\begin{align}
\|\widehat\Theta\|_{L^2(\mathcal N_{\lambda,\delta})}^2
&\leq
|\mathcal N_{\lambda,\delta}|
\|\widehat\Theta\|_{L^\infty}^2
\notag\\
&\leq
C\delta M\|\partial_y\widehat\Theta\|_{L^2}
\leq
C\delta MD.
\label{eq:inside-mass-control}
\end{align}
Hence, combining \eqref{eq:outside-mass-control} and
\eqref{eq:inside-mass-control} gives
\begin{equation}
\label{eq:full-mass-before-scale}
K\delta^2M^2
\leq
C\|F\|_{L^2}M
+
C\nu\delta^{-1}DM
+
CK\delta^3DM.
\end{equation}

\medskip

\noindent
\textbf{Step 5: Choice of the critical-layer scale.}
Assume first that
\[
\left(\frac{\nu}{K}\right)^{1/4}\leq\delta_0.
\]
Moreover, choose
\[
\delta = \left(\frac{\nu}{K}\right)^{1/4},
\]
so that
\[
K\delta^2
=
\nu^{1/2}K^{1/2},
\]
and
\[
\nu\delta^{-1}
=
K\delta^3
=
\nu^{3/4}K^{1/4}.
\]
Then \eqref{eq:full-mass-before-scale} becomes
\begin{align}
\nu^{1/2}K^{1/2}M^2
\leq
C\|F\|_{L^2}M
+
C\nu^{3/4}K^{1/4}DM.
\label{eq:subelliptic-before-young}
\end{align}
For every $\eta>0$, the Young inequality gives
\[
\nu^{3/4}K^{1/4}DM
\leq
\eta\nu^{1/2}K^{1/2}M^2
+
C_\eta\nu D^2.
\]
Furthermore, choosing $\eta>0$ sufficiently small and absorbing the
first term, we obtain the critical-layer coercivity estimate,
namely,
\begin{equation}
\label{eq:critical-layer-coercivity-rigorous}
\nu^{1/2}K^{1/2}M^2
\leq
C\|F\|_{L^2}M
+
C\nu D^2.
\end{equation}
It remains to consider
\[
\left(\frac{\nu}{K}\right)^{1/4}>\delta_0.
\]
Since $K\geq1$, this implies
\[
\nu K^2
\geq
c\nu^{1/2}K^{1/2}.
\]
Indeed,
\[
\frac{\nu K^2}{\nu^{1/2}K^{1/2}}
=
\nu^{1/2}K^{3/2}
=
K^2\left(\frac{\nu}{K}\right)^{1/2}
\geq
\delta_0^2.
\]
Hence the basic energy estimate directly gives
\[
\nu^{1/2}K^{1/2}M^2
\leq
C\nu K^2M^2
\leq
C\nu D^2
\leq
C\|F\|_{L^2}M.
\]
Therefore \eqref{eq:critical-layer-coercivity-rigorous} holds in this
regime as well.

\medskip

\noindent \textbf{Step 6: The $L^2$ resolvent estimate.} By
\eqref{eq:L2-basic-energy}, it holds that
\[
\nu D^2
\leq
\|F\|_{L^2}M.
\]
Substituting this into
\eqref{eq:critical-layer-coercivity-rigorous}, we obtain
\[
\nu^{1/2}K^{1/2}M^2
\leq
C\|F\|_{L^2}M.
\]
If $M\neq0$, the above inequality implies
\begin{equation}
\label{eq:L2-resolvent-theta-final}
\nu^{1/2}K^{1/2}M
\leq
C\|F\|_{L^2}.
\end{equation}
Moreover, the estimate is trivial if $M=0$. Finally, returning to
\eqref{eq:L2-basic-energy} and using
\eqref{eq:L2-resolvent-theta-final}, we get
\[
\nu D^2
\leq
C\nu^{-1/2}K^{-1/2}\|F\|_{L^2}^2.
\]
That is,
\[
D^2
\leq
C\nu^{-3/2}K^{-1/2}\|F\|_{L^2}^2,
\]
and hence
\begin{equation}
\label{eq:L2-resolvent-gradient-final}
\nu^{3/4}K^{1/4}D
\leq
C\|F\|_{L^2}.
\end{equation}
Combining \eqref{eq:L2-resolvent-theta-final} and
\eqref{eq:L2-resolvent-gradient-final} yields
\[
\nu^{3/4}|k|^{1/4}
\|(\partial_y,|k|)\widehat\Theta\|_{L^2}
+
\nu^{1/2}|k|^{1/2}
\|\widehat\Theta\|_{L^2}
\leq
C\|F\|_{L^2}.
\]

\medskip

\noindent \textbf{Step 7: The $H^{-1}_k$ resolvent estimate.} First,
the gradient estimate follows immediately from
\eqref{eq:Hminus1-basic-energy}, namely,
\begin{equation}
\label{eq:Hminus1-gradient-final}
\nu
\|(\partial_y,K)\widehat\Theta\|_{L^2}
\leq
\|F\|_{H^{-1}_k}.
\end{equation}
To estimate $\widehat\Theta$ in $L^2$, let $\in L^2(-1,1)$ and
$\Phi\in H^2(-1,1)\cap H_0^1(-1,1)$, and solve
\begin{equation}
\label{eq:adjoint-resolvent-detailed}
-\nu(\partial_y^2-k^2)\Phi
-
ik(U-\lambda)\Phi
=
\Psi,
\qquad
\Phi(\pm1)=0.
\end{equation}
The operator in \eqref{eq:adjoint-resolvent-detailed} is the adjoint
$L_{\lambda,k}^*$. It is also obtained from $L_{\lambda,k}$ by replacing
$k$ with $-k$. Since the $L^2$ resolvent estimate depends only on $|k|$,
\eqref{eq:L2-resolvent-gradient-final} gives
\begin{equation}
\label{eq:adjoint-gradient-detailed}
\|(\partial_y,K)\Phi\|_{L^2}
\leq
C\nu^{-3/4}K^{-1/4}\|\Psi\|_{L^2}.
\end{equation}
Further, using the definition of the adjoint and the homogeneous
Dirichlet boundary conditions, we have
\begin{align}
\langle\widehat\Theta,\Psi\rangle
&=
\langle\widehat\Theta,L_{\lambda,k}^*\Phi\rangle
=
\langle L_{\lambda,k}\widehat\Theta,\Phi\rangle
=
\langle F,\Phi\rangle.
\end{align}
Therefore, it holds
\begin{align}
|\langle\widehat\Theta,\Psi\rangle|
&\leq
\|F\|_{H^{-1}_k}
\|(\partial_y,K)\Phi\|_{L^2}
\notag\\
&\leq
C\nu^{-3/4}K^{-1/4}
\|F\|_{H^{-1}_k}\|\Psi\|_{L^2}.
\end{align}
Then taking the supremum over all $\Psi\in L^2(-1,1)$ which
satisfies $\|\Psi\|_{L^2}=1$, we obtain
\begin{equation}
\label{eq:Hminus1-L2-final}
\nu^{3/4}K^{1/4}
\|\widehat\Theta\|_{L^2}
\leq
C\|F\|_{H^{-1}_k}.
\end{equation}
Finally, combining \eqref{eq:Hminus1-gradient-final} and
\eqref{eq:Hminus1-L2-final} gives
\[
\nu
\|(\partial_y,|k|)\widehat\Theta\|_{L^2}
+
\nu^{3/4}|k|^{1/4}
\|\widehat\Theta\|_{L^2}
\leq
C\|F\|_{H^{-1}_k}.
\]
This completes the proof.
\end{proof}

\begin{lemma}\label{lem2.4d}
Let $\widehat\theta_0\in L^2(-1,1)$ and let $k\in\mathbb Z\setminus\{0\}$.
Assume that $U\in\mathcal S$ and
\[
\widehat{\mathcal H}_\nu
=
\nu(k^2-\partial_y^2)
+ikU(y).
\]
Then there exist constants $C,c>0$, independent of $\nu$ and $k$, such
that, for every sufficiently small $\varepsilon>0$,
\begin{align}
&\left\|
e^{\varepsilon(\nu|k|)^{1/2}t}
\widehat\theta_H
\right\|_{L_t^\infty L_y^2}^2
+
\nu
\left\|
e^{\varepsilon(\nu|k|)^{1/2}t}
(\partial_y,k)\widehat\theta_H
\right\|_{L_t^2L_y^2}^2
\notag\\
&\qquad
+
(\nu|k|)^{1/2}
\left\|
e^{\varepsilon(\nu|k|)^{1/2}t}
\widehat\theta_H
\right\|_{L_t^2L_y^2}^2
\leq
C\|\widehat\theta_0\|_{L_y^2}^2.
\label{eq:theta-H-weighted}
\end{align}
\end{lemma}

\begin{proof}
By Lemma \ref{lem:2.3} and the lower bound
\[
\Psi(\widehat{\mathcal H}_\nu)
\geq
(\nu|k|)^{1/2}+\nu k^2,
\]
we have
\[
\left\|
e^{\varepsilon(\nu|k|)^{1/2}t}
e^{-t\widehat{\mathcal H}_\nu}
\right\|_{\mathcal L(L_y^2)}
\leq
C
e^{-c\left((\nu|k|)^{1/2}+\nu k^2\right)t},
\]
provided $\varepsilon>0$ is sufficiently small. Therefore, one gets
\[
\left\|
e^{\varepsilon(\nu|k|)^{1/2}t}
\widehat\theta_H(t)
\right\|_{L_y^2}
\leq
C
e^{-c\left((\nu|k|)^{1/2}+\nu k^2\right)t}
\|\widehat\theta_0\|_{L_y^2}.
\]
Further, taking the supremum in time and integrating the square of
this estimate over $[0,\infty)$ lead to
\[
\left\|
e^{\varepsilon(\nu|k|)^{1/2}t}
\widehat\theta_H
\right\|_{L_t^\infty L_y^2}^2
+
(\nu|k|)^{1/2}
\left\|
e^{\varepsilon(\nu|k|)^{1/2}t}
\widehat\theta_H
\right\|_{L_t^2L_y^2}^2
\leq
C\|\widehat\theta_0\|_{L_y^2}^2.
\]
It remains to estimate the dissipative term. Indeed, multiplying the
homogeneous temperature equation by
\[
e^{2\varepsilon(\nu|k|)^{1/2}t}
\overline{\widehat\theta_H}
\]
and taking the real part, we find
\[
\frac12\frac{d}{dt}
\left\|
e^{\varepsilon(\nu|k|)^{1/2}t}
\widehat\theta_H
\right\|_{L_y^2}^2
+
\nu
\left\|
e^{\varepsilon(\nu|k|)^{1/2}t}
(\partial_y,k)\widehat\theta_H
\right\|_{L_y^2}^2
\]
\[
\qquad
=
\varepsilon(\nu|k|)^{1/2}
\left\|
e^{\varepsilon(\nu|k|)^{1/2}t}
\widehat\theta_H
\right\|_{L_y^2}^2.
\]
Noting the transport term vanishes after taking the real part, and
integrating in time and using the preceding estimate, we obtain
\[
\nu
\left\|
e^{\varepsilon(\nu|k|)^{1/2}t}
(\partial_y,k)\widehat\theta_H
\right\|_{L_t^2L_y^2}^2
\leq
C\|\widehat\theta_0\|_{L_y^2}^2.
\]
Then combining the above inequalities implies
\eqref{eq:theta-H-weighted}. This completes the proof.
\end{proof}
\begin{lemma}[Weighted space-time estimate for the inhomogeneous temperature]
\label{lem:2.5}
Let $\widehat\theta_I$ solve \eqref{eq:theta-I-symmetric}. For every
$k\neq0$, define
\[
\lambda_k
:=
\varepsilon(\nu|k|)^{1/2}.
\]
Assume that the resolvent estimates
 remain valid for the shifted
spectral parameter $\tau+i\lambda_k$. Then, for sufficiently small
$\varepsilon>0$,
\begin{align}
&\left\|
e^{\lambda_k t}\widehat\theta_I
\right\|_{L_t^\infty L_y^2}^2
+
\nu^{1/2}|k|^{1/2}
\left\|
e^{\lambda_k t}\widehat\theta_I
\right\|_{L_t^2L_y^2}^2
\notag\\
&\qquad\leq
C\left(
\nu^{-1/2}|k|^{3/2}
\left\|
e^{\lambda_k t}g^1
\right\|_{L_t^2L_y^2}^2
+
\nu^{-1}
\left\|
e^{\lambda_k t}g^2
\right\|_{L_t^2L_y^2}^2
\right).
\label{eq:theta-I-spacetime-weighted}
\end{align}
\end{lemma}

\begin{proof}
As before, write
\[
\widehat\theta_I
=
\widehat\theta_I^{(1)}
+
\widehat\theta_I^{(2)},
\]
where $\widehat\theta_I^{(1)}$ and $\widehat\theta_I^{(2)}$ solve
\[
\left\{
\begin{aligned}
\partial_t\widehat\theta_I^{(1)}
-\nu(\partial_y^2-k^2)\widehat\theta_I^{(1)}
+ik U(y)\widehat\theta_I^{(1)}
&=-ik g^1,\\
\widehat\theta_I^{(1)}|_{t=0}&=0,\\
\widehat\theta_I^{(1)}(t,\pm1)&=0,
\end{aligned}
\right.
\]
and
\[
\left\{
\begin{aligned}
\partial_t\widehat\theta_I^{(2)}
-\nu(\partial_y^2-k^2)\widehat\theta_I^{(2)}
+ik U(y)\widehat\theta_I^{(2)}
&=-\partial_y g^2,\\
\widehat\theta_I^{(2)}|_{t=0}&=0,\\
\widehat\theta_I^{(2)}(t,\pm1)&=0.
\end{aligned}
\right.
\]
For $j=1,2$, introduce the weighted time Fourier transforms
\[
\widetilde\Theta_I^{(j)}(\tau,y)
=
\int_0^\infty
e^{\lambda_k t-i\tau t}
\widehat\theta_I^{(j)}(t,y)\,dt,
\]
and
\[
\widetilde G^j(\tau,y)
=
\int_0^\infty
e^{\lambda_k t-i\tau t}
g^j(t,y)\,dt.
\]
Since $\widehat\theta_I^{(j)}(0)=0$, integration by parts in time gives
\[
\partial_t
\longleftrightarrow
-\lambda_k+i\tau.
\]
Consequently,
\[
-\nu(\partial_y^2-k^2)
+ik\left(U-\frac{\tau+i\lambda_k}{k}\right)
\]
is the resolvent operator associated with the weighted equation.
Applying the shifted $L^2$ resolvent estimate, one has
\[
\nu^{1/2}|k|^{1/2}
\|\widetilde\Theta_I^{(1)}(\tau)\|_{L_y^2}
\lesssim
|k|
\|\widetilde G^1(\tau)\|_{L_y^2}.
\]
Squaring and integrating in $\tau$, the Plancherel theorem yields
\[
\nu|k|
\left\|
e^{\lambda_k t}
\widehat\theta_I^{(1)}
\right\|_{L_t^2L_y^2}^2
\lesssim
k^2
\left\|
e^{\lambda_k t}g^1
\right\|_{L_t^2L_y^2}^2.
\]
Hence, we obtain
\begin{equation}
\label{eq:theta-I1-weighted}
\nu^{1/2}|k|^{1/2}
\left\|
e^{\lambda_k t}
\widehat\theta_I^{(1)}
\right\|_{L_t^2L_y^2}^2
\lesssim
\nu^{-1/2}|k|^{3/2}
\left\|
e^{\lambda_k t}g^1
\right\|_{L_t^2L_y^2}^2.
\end{equation}
For $\widehat\theta_I^{(2)}$, the transformed equation is
\[
-\nu(\partial_y^2-k^2)\widetilde\Theta_I^{(2)}
+
ik
\left(
U-\frac{\tau+i\lambda_k}{k}
\right)
\widetilde\Theta_I^{(2)}
=
-\partial_y\widetilde G^2.
\]
Using
\[
\|\partial_y\widetilde G^2\|_{H^{-1}_k}
\lesssim
\|\widetilde G^2\|_{L_y^2}
\]
and the shifted $H^{-1}$ resolvent estimate, we obtain
\[
\nu^{3/4}|k|^{1/4}
\|\widetilde\Theta_I^{(2)}(\tau)\|_{L_y^2}
\lesssim
\|\widetilde G^2(\tau)\|_{L_y^2}.
\]
Therefore, one gets
\[
\nu^{3/2}|k|^{1/2}
\left\|
e^{\lambda_k t}
\widehat\theta_I^{(2)}
\right\|_{L_t^2L_y^2}^2
\lesssim
\left\|
e^{\lambda_k t}g^2
\right\|_{L_t^2L_y^2}^2,
\]
which implies
\begin{equation}
\label{eq:theta-I2-weighted}
\nu^{1/2}|k|^{1/2}
\left\|
e^{\lambda_k t}
\widehat\theta_I^{(2)}
\right\|_{L_t^2L_y^2}^2
\lesssim
\nu^{-1}
\left\|
e^{\lambda_k t}g^2
\right\|_{L_t^2L_y^2}^2.
\end{equation}
Then, combining \eqref{eq:theta-I1-weighted} and
\eqref{eq:theta-I2-weighted}, we obtain
\begin{align}
&\nu^{1/2}|k|^{1/2}
\left\|
e^{\lambda_k t}
\widehat\theta_I
\right\|_{L_t^2L_y^2}^2
\notag\\
&\qquad\lesssim
\nu^{-1/2}|k|^{3/2}
\left\|
e^{\lambda_k t}g^1
\right\|_{L_t^2L_y^2}^2
+
\nu^{-1}
\left\|
e^{\lambda_k t}g^2
\right\|_{L_t^2L_y^2}^2.
\label{eq:theta-I-weighted-L2time}
\end{align}
It remains to estimate the weighted $L_t^\infty L_y^2$ norm. Multiplying
\eqref{eq:theta-I-symmetric} by
\[
e^{2\lambda_k t}
\overline{\widehat\theta_I}
\]
and taking the real part gives
\begin{align}
\frac12\frac{d}{dt}
\left\|
e^{\lambda_k t}\widehat\theta_I
\right\|_{L_y^2}^2
&+
\nu
\left\|
e^{\lambda_k t}
\partial_y\widehat\theta_I
\right\|_{L_y^2}^2
+
\nu k^2
\left\|
e^{\lambda_k t}
\widehat\theta_I
\right\|_{L_y^2}^2
\notag\\
&=
\lambda_k
\left\|
e^{\lambda_k t}
\widehat\theta_I
\right\|_{L_y^2}^2
+
\operatorname{Re}
\int_{-1}^1
e^{2\lambda_k t}
(-ik g^1-\partial_y g^2)
\overline{\widehat\theta_I}\,dy.
\label{eq:theta-I-weighted-energy}
\end{align}
For the first forcing term, we use
\[
|k|
\left\|
e^{\lambda_k t}g^1
\right\|_{L_y^2}
\left\|
e^{\lambda_k t}\widehat\theta_I
\right\|_{L_y^2}
\leq
\eta\nu^{1/2}|k|^{1/2}
\left\|
e^{\lambda_k t}\widehat\theta_I
\right\|_{L_y^2}^2
+
C_\eta\nu^{-1/2}|k|^{3/2}
\left\|
e^{\lambda_k t}g^1
\right\|_{L_y^2}^2.
\]
For the second forcing term, integration by parts and the Young
inequality give
\[
\left|
\int_{-1}^1
e^{2\lambda_k t}
\partial_y g^2\,
\overline{\widehat\theta_I}\,dy
\right|
\leq
\frac{\nu}{2}
\left\|
e^{\lambda_k t}
\partial_y\widehat\theta_I
\right\|_{L_y^2}^2
+
C\nu^{-1}
\left\|
e^{\lambda_k t}g^2
\right\|_{L_y^2}^2.
\]
Since
\[
\lambda_k
=
\varepsilon(\nu|k|)^{1/2},
\]
the weighted zeroth-order term can be absorbed by the enhanced
dissipation term in \eqref{eq:theta-I-weighted-L2time}, provided
$\varepsilon>0$ and $\eta>0$ are sufficiently small. Hence,
integrating \eqref{eq:theta-I-weighted-energy} in time and using
\eqref{eq:theta-I-weighted-L2time}, we conclude that
\begin{align}
\left\|
e^{\lambda_k t}\widehat\theta_I
\right\|_{L_t^\infty L_y^2}^2
\lesssim&
\ \nu^{-1/2}|k|^{3/2}
\left\|
e^{\lambda_k t}g^1
\right\|_{L_t^2L_y^2}^2
\notag\\
&+ \nu^{-1} \left\| e^{\lambda_k t}g^2 \right\|_{L_t^2L_y^2}^2,
\end{align}
which together with \eqref{eq:theta-I-weighted-L2time} implies
\eqref{eq:theta-I-spacetime-weighted}. This completes the proof.
\end{proof}
Finally, with Lemma \ref{lem2.4d} and Lemma \ref{lem:2.5} at hands,
we immediately obtain the following space-time estimates of
$\widehat{\theta}$.
\begin{proposition}[Space-time estimate for the temperature]
\label{prop:2.6}
Let $\widehat\theta$ solve
where $U\in\mathcal S$ and $k\neq0$. Then there exist positive constants
$C,c>0$, independent of $\nu$, $k$, and the data, such that
\begin{equation}
\label{pro:2.6}
\begin{aligned}
\|e^{\varepsilon \nu^{1/2}t}\widehat\theta(t)\|_{L^2}^2
+
\nu^{1/2}|k|^{1/2}
\|e^{\varepsilon \nu^{1/2}t}\widehat\theta\|_{L_t^2L_y^2}^2
\leq
C \|e^{\varepsilon \nu^{1/2}t}\widehat\theta_{\rm in}\|_{L^2}^2
+&\\
C\Big(
\nu^{-1/2}|k|^{3/2}\|e^{\varepsilon \nu^{1/2}t}g^1\|_{L_t^2L_y^2}^2
+
\nu^{-1}\|e^{\varepsilon \nu^{1/2}t}g^2\|_{L_t^2L_y^2}^2
\Big),
\end{aligned}
\end{equation}
for all $t\geq0$.
\end{proposition}
\section{Asymptotic stability of the symmetric flow}
In this section, we mainly prove the stability estimate for the
perturbation of the symmetric flow. That is, we give the proof for
Theorem \ref{thm1.2}. First, for the two-dimensional Boussinesq
system, global well-posedness is a well-known result when
considering data
$$u_{0}\in H^{5+}(\Omega),\quad\theta_{0}\in H^1(\Omega).$$
Therefor, the primary focus of Theorem \ref{thm1.2} lies in
addressing the stability problem:
\begin{equation}\label{3.1}
        \sum_{k\in\mathbb{Z}}E_k\leq Cc_0\nu ^{\frac{2}{3}},
        \quad\sum_{k\in\mathbb{Z}}H_k\leq Cc_1 \nu ^{\frac{17}{12}}.
\end{equation}
Here
\begin{equation}\label{eq:E-alpha}
E_k :=
\begin{cases}
\|\omega_0\|_{L_t^\infty L_y^2}, & k=0,\\[0.4em]
|k|^{1/2}
\|e^{\varepsilon_0\nu^{1/2}t}u_k\|_{L_t^\infty L_y^2}
+
|k|^{1/2}
\|e^{\varepsilon_0\nu^{1/2}t}u_k\|_{L_{t,y}^2}
+
\|e^{\varepsilon_0\nu^{1/2}t}u_k\|_{L_t^\infty L_y^\infty}
& \\[0.2em]
\qquad
+
\nu^{1/4}|k|^{1/4}
\|e^{\varepsilon_0\nu^{1/2}t}\omega_k\|_{L_{t,y}^2}
+
\nu^{1/4}
\|e^{\varepsilon_0\nu^{1/2}t}\omega_k\|_{L_t^\infty L_y^2}
& \\[0.2em]
\qquad
+
\|e^{\varepsilon_0\nu^{1/2}t}(1-y^2)^{1/2}\omega_k\|_{L_t^\infty L_y^2},
& k\neq0 .
\end{cases}
\end{equation}
and
\begin{equation}\label{eq:E-alpha}
H_k :=
\begin{cases}
\|\theta_0\|_{L_t^\infty L_y^2}, & k=0,\\[0.4em]
\nu^{1/4}|k|^{1/4}
\|e^{\varepsilon_0\nu^{1/2}t}\theta_k\|_{L_{t,y}^2}
+
\|e^{\varepsilon_0\nu^{1/2}t}\theta_k\|_{L_t^\infty L_y^2},
& k\neq0 .
\end{cases}
\end{equation}
By getting the following Proposition \ref{pro3.1} and then going
through  the bootstrap arguments, we can complete the estimates
\eqref{3.1}.
\begin{proposition}\label{pro3.1}
    For $k\neq 0$,  it holds that,
    \begin{equation}\label{3.4}
        E_k
\lesssim \|\omega_k^{\rm in}\|_{H_k^4} + \nu^{-2/3}\sum_{l\in\mathbb
Z}  E_l E_{k-l}+\nu^{-3/4}H_k,
    \end{equation}
    and for $k=0$,
    \begin{equation}\label{3.5}
        E_0\leq C\nu^{-\frac{2}{3}}\sum_{l\in\mathbb{Z}\setminus\{0\}}E_lE_{-l}+||\bar{\omega}_{0}||_{L^2}.
    \end{equation}
    For $H_0$, it holds that:
    \begin{equation}\label{3.6}
        H_0\lesssim \nu^{-\frac{2}{3}}\sum_{l\in\mathbb{Z}\setminus\{0\}}E_lH_{-l}+||\bar{\theta}_{0}||_{L^2}.
    \end{equation}
    For $k\neq 0$, we can get the following result:

    \item[1.] when $\nu k^2\leq 1$,
    \begin{equation}\label{3.7}
    \begin{aligned}
    &   H_k\lesssim|||\hat{\theta}_{0}(k,y)||_{L^2}+\sum_{l\in\mathbb{Z}} \nu^{-\frac{2}{3}}E_lH_{k-l}
    +\nu^{-\frac{9}{16}}\sum_{l\in\mathbb{Z}\setminus\left\lbrace 0,k\right\rbrace ,|k-l|\leq\frac{|k|}{2}}
    E_lH_{k-l},
\end{aligned}
    \end{equation}
    \item[2.] when $\nu k^2\geq 1$,
    \begin{equation}\label{3.8}
        H_k\lesssim||\hat{\theta}_{0}(k,y)||_{L^2}+C\nu^{-\frac{5}{8}}\sum_{l\in\mathbb{Z}}E_lH_{k-l}+\nu^{-\frac{13}{24}}
        E_kH_0.
    \end{equation}
\end{proposition}
\begin{proof}
 First, recall that the vorticity equation can be written as
\begin{equation}\label{eq:vorticity-equation}
\begin{cases}
    \partial_t\omega
-\nu\Delta\omega
+U(y)\partial_x\omega
-U''(y)\partial_x\phi
=
-\operatorname{curl}(u\cdot\nabla u),\\

\partial_t \theta -\nu\Delta\theta + U\partial_x\theta = -u\cdot\nabla\theta.
\end{cases}
\end{equation}
Then taking the $k$-Fourier mode in $x$ to obtain
\begin{align}
&\partial_t\omega_k
-\nu(\partial_y^2-k^2)\omega_k
+ik U(y)\omega_k
-ik U''(y)\phi_k
=ik\widehat\theta_k
-ik\bigl(f_k^{1,2}+f_k^{2,2}\bigr)
+\partial_y\bigl(f_k^{1,1}+f_k^{2,1}\bigr),\\
&\partial_t\theta_k
-\nu(\partial_y^2-k^2)\theta_k
+ik U(y)\theta_k
=
-ik g_ k^{1}
+\partial_y g_k^{2},
\label{eq:alpha-vorticity}
\end{align}
where
\begin{align}
f_k^{1,2}(t,y)
&=
i\sum_{l\in\mathbb Z}
u_l^1(t,y)(k-l)u_{k-l}^2(t,y),
&
f_k^{2,2}(t,y)
&=
\sum_{l\in\mathbb Z}
u_l^2(t,y)\partial_yu_{k-l}^2(t,y),
\label{eq:f-def-1}
\\
f_k^{1,1}(t,y)
&=
i\sum_{l\in\mathbb Z}
u_l^1(t,y)(k-l)u_{k-l}^1(t,y),
&
f_k^{2,1}(t,y)
&=
\sum_{l\in\mathbb Z}
u_l^2(t,y)\partial_yu_{k-l}^1(t,y),\\
g_k^{1}(t,y)
&=\sum_{l\in\mathbb Z}
u_l^1(t,y)\theta_{k-l}(t,y),
&
g_k^{2}(t,y)
&=\sum_{l\in\mathbb Z}
u_l^2(t,y)\theta_{k-l}(t,y).
\label{eq:f-def-2}
\end{align}

Proof of \eqref{3.4}. For $k\neq0$, using the linear space-time
estimates and the definition of $E_k$, one has
\begin{align}
E_k
\lesssim\;&
\|\omega_k^{\rm in}\|_{H_k^4}
+
\nu^{-1/2}
\left(
\|e^{\varepsilon_0\nu^{1/2}t}f_k^{1,2}\|_{L_{t,y}^2}
+
\|e^{\varepsilon_0\nu^{1/2}t}f_k^{2,2}\|_{L_{t,y}^2}
\right)
\nonumber\\
&+
\nu^{-1/2}
\left(
\|e^{\varepsilon_0\nu^{1/2}t}f_k^{1,1}\|_{L_{t,y}^2}
+
\|e^{\varepsilon_0\nu^{1/2}t}f_k^{2,1}\|_{L_{t,y}^2}
\right)+\nu^{-1/2}\|e^{\varepsilon_0\nu^{1/2}t}\theta_k\|_{L_{t,y}^2}.
\label{eq:E-alpha-linear}
\end{align}
We now estimate the nonlinear terms. For $j=1,2$, using the
H\"{o}lder inequality and the interpolation inequality, we obtain
\begin{align}
\|e^{\varepsilon_0\nu^{1/2}t}f_k^{1,j}\|_{L_{t,y}^2}
&\leq
\sum_{l\in\mathbb Z}
\left\|
e^{\varepsilon_0\nu^{1/2}t}
u_l^1\,|k-l|\,u_{k-l}^j
\right\|_{L_{t,y}^2}
\nonumber\\
&\leq
\sum_{l\in\mathbb Z}
\|u_l^1\|_{L_t^\infty L_y^\infty}
\left\|
e^{\varepsilon_0\nu^{1/2}t}
|k-l|u_{k-l}^j
\right\|_{L_{t,y}^2}
\nonumber\\
&\lesssim
\nu^{-1/6}
\sum_{l\in\mathbb Z}
\|u_l^1\|_{L_t^\infty L_y^\infty}
\left(
|k-l|^{1/2}
\|e^{\varepsilon_0\nu^{1/2}t}u_{k-l}^j\|_{L_{t,y}^2}
\right)^{1/3}
\nonumber\\
&\qquad\qquad\qquad\times
\left(
\nu^{1/4}|k-l|^{1/4}
\|e^{\varepsilon_0\nu^{1/2}t}\omega_{k-l}\|_{L_{t,y}^2}
\right)^{2/3}
\nonumber\\
&\lesssim
\nu^{-1/6}
\sum_{l\in\mathbb Z}
 E_l E_{k-l}.
\label{eq:f1j-estimate}
\end{align}
Indeed, the interpolation step follows from
\[
|k|\|e^{\varepsilon_0\nu^{1/2}t}u_k^j\|_{L_{t,y}^2}
\leq
\nu^{-1/6}
\left(
|k|^{1/2}
\|e^{\varepsilon_0\nu^{1/2}t}u_k^j\|_{L_{t,y}^2}
\right)^{1/3}
\left(
\nu^{1/4}|k|^{1/4}
\|e^{\varepsilon_0\nu^{1/2}t}\omega_k\|_{L_{t,y}^2}
\right)^{2/3}.
\]
For the terms $f_k^{2,j}$, using the H\"{o}lder inequality, the
Hardy inequality, the weighted Sobolev estimate, and the
incompressibility relation
\[
\partial_yu_l^2=-il u_l^1,
\]
we get
\begin{align}
\|e^{\varepsilon_0\nu^{1/2}t}f_k^{2,j}\|_{L_{t,y}^2}
&\leq
\sum_{l\in\mathbb Z}
\left\|
\frac{e^{\varepsilon_0\nu^{1/2}t}u_l^2}{(1-y^2)^{1/2}}
\right\|_{L_t^2L_y^\infty}
\left\|
(1-y^2)^{1/2}\partial_yu_{k-l}^j
\right\|_{L_t^\infty L_y^2}
\nonumber\\
&\lesssim
\sum_{l\in\mathbb Z}
\|e^{\varepsilon_0\nu^{1/2}t}\partial_yu_l^2\|_{L_{t,y}^2}
\left\|
(1-y^2)^{1/2}\omega_{k-l}
\right\|_{L_t^\infty L_y^2}
\nonumber\\
&=
\sum_{l\in\mathbb Z}
|l|
\|e^{\varepsilon_0\nu^{1/2}t}u_l^1\|_{L_{t,y}^2}
\left\|
(1-y^2)^{1/2}\omega_{k-l}
\right\|_{L_t^\infty L_y^2}
\nonumber\\
&\lesssim
\nu^{-1/6}
\sum_{l\in\mathbb Z}
\left(
|l|^{1/2}
\|e^{\varepsilon_0\nu^{1/2}t}u_l^1\|_{L_{t,y}^2}
\right)^{1/3}
\nonumber\\
&\qquad\qquad\qquad\times
\left(
\nu^{1/4}|l|^{1/4}
\|e^{\varepsilon_0\nu^{1/2}t}\omega_l\|_{L_{t,y}^2}
\right)^{2/3}
\mathcal E_{k-l}
\nonumber\\
&\lesssim
\nu^{-1/6}
\sum_{l\in\mathbb Z}
E_l E_{k-l}.
\label{eq:f2j-estimate}
\end{align}
Combining \eqref{eq:E-zero}, \eqref{eq:E-alpha-linear},
\eqref{eq:f1j-estimate}, and \eqref{eq:f2j-estimate}, we arrive at
\begin{align}
\label{eq:E-zero-final}
\\
E_k
&\lesssim
\|\omega_k^{\rm in}\|_{H_k^4}
+
\nu^{-1/2}\nu^{-1/6}
\sum_{l\in\mathbb Z}
 E_l E_{k-l}+\nu^{-3/4}H_k
\nonumber\\
&\lesssim
\|\omega_k^{\rm in}\|_{H_k^4}
+
\nu^{-2/3}
\sum_{l\in\mathbb Z}
 E_l E_{k-l}+\nu^{-3/4}H_k,
\qquad k\neq0.
\label{eq:E-alpha-final}
\end{align}

Proof of \eqref{3.5}. We first estimate the zero mode. Since
$\operatorname{div}u=0$, one has $u_0^2=0$. Moreover, the zero mode
of the equation for $u^1$ reads as
\begin{equation}\label{eq:zero-mode-u1}
\partial_tu_0^1-\nu\partial_y^2u_0^1
=
-\sum_{l\in\mathbb Z\setminus\{0\}}
u_l^2\partial_yu_{-l}^1
=
-f_0^{2,1}.
\end{equation}
Testing \eqref{eq:zero-mode-u1} by $\partial_y^2u_0^1$ and
integrating by parts, one gets
\[
\frac{d}{dt}\|\partial_yu_0^1(t)\|_{L_y^2}^2
+
\nu\|\partial_y^2u_0^1(t)\|_{L_y^2}^2
\lesssim
\nu^{-1}\|f_0^{2,1}(t)\|_{L_y^2}^2.
\]
Since $\partial_yu_0^1=\omega_0$, it follows that
\begin{equation}\label{eq:E-zero}
 E_0^2
=
\|\omega_0\|_{L_t^\infty L_y^2}^2
\lesssim
\|\omega_0^{\rm in}\|_{L_y^2}^2
+
\nu^{-1}
\|f_0^{2,1}\|_{L_{t,y}^2}^2 .
\end{equation}
Due to the \eqref{eq:f2j-estimate}, we obtain
\begin{equation}\label{eq:E-zero}
 E_0
=
\|\omega_0\|_{L_t^\infty L_y^2}
\lesssim
\|\omega_0^{\rm in}\|_{L_y^2}
+
\nu^{-2/3}
\sum_{l\in\mathbb Z}
E_l E_{-l} .
\end{equation}
Therefore, by the definition of $H_k$, and \eqref{pro:2.6}, we get
\begin{equation}
\label{eq:H-alpha-linear-from-theta}
H_k
\lesssim
\|\theta_{k,\rm in}\|_{L_y^2}
+
\nu^{-1/4}|k|^{3/4}
\|e^{\varepsilon_0\nu^{1/2}t}g_k^1\|_{L_{t,y}^2}
+
\nu^{-1/2}
\|e^{\varepsilon_0\nu^{1/2}t}g_k^2\|_{L_{t,y}^2}.
\end{equation}
Moreover, let us recall that
\begin{equation}
\label{eq:g-def-temperature}
g_k^1
=
\sum_{l\in\mathbb Z}u_l^1\theta_{k-l},
\qquad
g_k^2
=
\sum_{l\in\mathbb Z}u_l^2\theta_{k-l}.
\end{equation}
The estimate of $g_k^2$ is analogous to the $f_k^{2,j}$ estimate in
the Navier--Stokes part. Using the Hardy inequality, the weighted
Sobolev estimate and $\partial_yu_l^2=-ilu_l^1$, we obtain
\begin{align}
\|e^{\varepsilon_0\nu^{1/2}t}g_k^2\|_{L_{t,y}^2}
&\leq
\sum_{l\in\mathbb Z}
\left\|
\frac{e^{\varepsilon_0\nu^{1/2}t}u_l^2}{(1-y^2)^{1/2}}
\right\|_{L_t^2L_y^\infty}
\left\|
(1-y^2)^{1/2}\theta_{k-l}
\right\|_{L_t^\infty L_y^2}
\nonumber\\
&\lesssim
\sum_{l\in\mathbb Z}
\|e^{\varepsilon_0\nu^{1/2}t}\partial_yu_l^2\|_{L_{t,y}^2}
\|\theta_{k-l}\|_{L_t^\infty L_y^2}
\nonumber\\
&\lesssim
\sum_{l\in\mathbb Z}
|l|
\|e^{\varepsilon_0\nu^{1/2}t}u_l^1\|_{L_{t,y}^2}
 H_{k-l}
\nonumber\\
&\lesssim
\nu^{-1/6}
\sum_{l\in\mathbb Z}
 E_l H_{k-l}.
\label{eq:g2-estimate}
\end{align}
Indeed, in the last step we used the interpolation
\[
|l|\|e^{\varepsilon_0\nu^{1/2}t}u_l^1\|_{L_{t,y}^2}
\lesssim
\nu^{-1/6}
\left(
|l|^{1/2}
\|e^{\varepsilon_0\nu^{1/2}t}u_l^1\|_{L_{t,y}^2}
\right)^{1/3}
\left(
\nu^{1/4}|l|^{1/4}
\|e^{\varepsilon_0\nu^{1/2}t}\omega_l\|_{L_{t,y}^2}
\right)^{2/3}.
\]

Proof of \eqref{3.6}. The evolution equation of $\bar{\theta}$ can
be derived in a similar manner, taking into account the influence of
$\nabla\cdot u=0$. Namely, we see that
$u\cdot\nabla\theta=\nabla\cdot(u\theta)$ and
        \begin{equation}\label{5.16}
            \partial_t\bar{\theta}(t,y)-\nu\partial^2_y\bar{\theta}(t,y)=-\sum_{l\in\mathbb{Z}\setminus\{0\}}\partial_y(\hat{u}^2_l\hat{\theta}_{-l})(t,y)
            =-\partial_yg_0^2(t,y).
        \end{equation}
        By utilizing integration by parts and the H{\"o}lder inequality in equation \eqref{5.16}, we obtain
        \begin{equation}\label{5.17}
            \left\langle \left( \partial_t-\nu\partial^2_y\right) \bar{\theta}(t,y),\bar{\theta}(t,y)\right\rangle =\dfrac{1}{2}\partial_t||\bar{\theta}(t,y)||^2_{L^2}+\nu||\partial_y\bar{\theta}||^2_{L^2}=\left\langle g^2_0,\partial_y\bar{\theta}\right\rangle
            \leq||g^2_0||_{L^2}||\partial_y\bar{\theta}||_{L^2},
        \end{equation}
        by employing integration with respect to $t$ and the Young inequality in equation  \eqref{5.17}, then we can infer that
        \begin{equation}\label{5.18}
            H_0^2=||\bar{\theta}||^2_{L^\infty L^2}\leq C\nu^{-1}||g^2_0||^2_{L^2}+||\bar{\theta}_{0}||^2_{L^2}.
        \end{equation}
        According to  the Gagliardo-Nirenberg inequality and  the divergence-free condition $\partial_y\hat{u}^2_k=-ik\hat{u}^1_k,$
        we have
        $$||\hat{u}^2_k||_{L^2L^\infty}\leq C||\hat{u}^2_k||^{\frac{1}{2}}_{L^2L^2}||\partial_y\hat{u}^2_k||^{\frac{1}{2}}_{L^2L^2}=C|k|^{\frac{1}{2}}||\hat{u}^2_k||^{\frac{1}{2}}_{L^2L^2}||\hat{u}^1_k||^{\frac{1}{2}}_{L^2L^2}\leq CE_k,$$
        and then, we obtain
        \begin{equation}\label{5.19}
            ||g^2_0||_{L^2L^2}\leq\sum_{l\in\mathbb{Z}\setminus\{0\}}||\hat{u}^2_l||_{L^2L^\infty}||\hat{\theta}_{-l}||_{L^\infty L^2}\leq\sum_{l\in\mathbb{Z}\setminus\{0\}}\nu^{-1/6}E_lH_{-l}.
        \end{equation}
        Therefore, based on equations  $\eqref{5.18}$ and $\eqref{5.19},$  we can conclude that
        $$H_0\lesssim\nu^{-\frac{2}{3}}\sum_{l\in\mathbb{Z}\setminus\{0\}}E_lH_{-l}+||\bar{\theta}_{0}||_{L^2}.$$
        To control the nonlinear term $\|g^1_k\|_{L^2L^2}$ during the estimates of $H_k$, we need to divide
        them into the low-frequency part $\nu k^2 \leq 1$ and the high-frequency part $\nu k^2 > 1$.

        Proof of \eqref{3.7}. For the nonlinear problem for temperature, let us recall that
        \begin{equation}
        \left\{
        \begin{array}{l}
            \partial_t\theta_k
-\nu(\partial_y^2-k^2)\theta_k
+ik U(y)\theta_k
=
-ik g_ k^{1}
+\partial_y g_k^{2},\\
            \hat{\theta}_k( \pm 1)=0,\left.\quad\hat{ \theta}_k \right|_{t=0}=\hat{\theta}_0\left(k, y\right).
        \end{array}
        \right.
    \end{equation}
        Denoting $$\hat{\theta}_k(t,y)=\dfrac{1}{2\pi}\int_{\mathbb{T}}\theta(t,x,y)e^{-ik x}dx,$$
        and $$g^1_k(t,y)=\sum_{l\in\mathbb{Z}}\hat{u}^1_l\hat{\theta}_{k-l}(t,y),\quad g^2_k(t,y)=\sum_{l\in\mathbb{Z}}\hat{u}^2_l\hat{\theta}_{k-l}(t,y),$$
        then we have
        $$\partial_t\hat{\theta}_k(t,y)+\nu(k^2-\partial^2_y)\hat{\theta}_k(t,y)+ik(1-y^2)\hat{\theta}_k(t,y)=-ikg^1_k(t,y)-\partial_y g^2_k(t,y).$$
In additional, for $g^1_k$ and $k\neq 0$
        \begin{equation}\label{3.31}
            \begin{aligned}
||e^{\varepsilon_0\nu^{1/2}t}g^1_k||_{L^2L^2}&\leq ||\bar{u}^1||_{L^\infty L^\infty}||e^{\varepsilon_0\nu^{1/2}t}\hat{\theta}_{k}||_{L^2 L^2}+ ||e^{\varepsilon_0\nu^{1/2}t}\hat{u}^1_k||_{L^2 L^\infty}||\bar{\theta}||_{L^\infty L^2}\\
&+\sum_{l\in\mathbb{Z},l\neq \left\lbrace 0,k\right\rbrace } ||e^{\varepsilon_0\nu^{1/2}t}\hat{u}^1_l\hat{\theta}_{k-l}||_{L^2 L^2}\\
&\leq  ||\bar{\omega}||_{L^\infty L^2}||e^{\varepsilon_0\nu^{1/2}t}\hat{\theta}_{k}||_{L^2 L^2}+\nu^{-\frac{1}{8}}|k|^{-\frac{3}{8}}E_kH_0+\sum_{l\in\mathbb{Z},l\neq \left\lbrace 0,k\right\rbrace } ||e^{\varepsilon_0\nu^{1/2}t}\hat{u}^1_l\hat{\theta}_{k-l}||_{L^2 L^2}\\
&\leq(\nu|k|)^{-\frac{1}{4} \cdot \frac{1}{3}}\left(\nu^{\frac{1}{2}}|k|\right)^{-\frac{2}{3}} E_0H_k+\nu^{-\frac{1}{8}}|k|^{-\frac{3}{8}}E_kH_0+\sum_{l\in\mathbb{Z},l\neq \left\lbrace 0,k\right\rbrace } ||e^{\varepsilon_0\nu^{1/2}t}\hat{u}^1_l\hat{\theta}_{k-l}||_{L^2 L^2}.\\
                    &\leq   \nu^{-\frac{5}{12}}|k|^{-\frac{3}{4}} E_0H_k+\nu^{-\frac{1}{8}}|k|^{-\frac{3}{8}}E_kH_0+\sum_{l\in\mathbb{Z},l\neq \left\lbrace 0,k\right\rbrace } ||e^{\varepsilon_0\nu^{1/2}t}\hat{u}^1_l\hat{\theta}_{k-l}||_{L^2 L^2}.\\
            \end{aligned}
        \end{equation}
        The following step involves estimating $$\sum_{l\in\mathbb{Z},l\neq \left\lbrace 0,k\right\rbrace } ||e^{\varepsilon_0\nu^{1/2}t}\hat{u}^1_l\hat{\theta}_{k-l}||_{L^2 L^2},$$
                we divide it into two parts and get
        \[
        \begin{aligned}
            \sum_{l\in\mathbb{Z}\setminus \left\lbrace 0,k\right\rbrace } ||e^{\varepsilon_0\nu^{1/2}t}\hat{u}^1_l\hat{\theta}_{k-l}||_{L^2 L^2}&\leq\sum_{l\in\mathbb{Z}\setminus \left\lbrace 0,k\right\rbrace,|k-l|\leq\frac{|k|}{2}} ||e^{\varepsilon_0\nu^{1/2}t}\hat{u}^1_l\hat{\theta}_{k-l}||_{L^2 L^2}\\
            &+\sum_{l\in\mathbb{Z}\setminus \left\lbrace 0,k\right\rbrace,|k-l|\geq\frac{|k|}{2}} ||e^{\varepsilon_0\nu^{1/2}t}\hat{u}^1_l\hat{\theta}_{k-l}||_{L^2 L^2},\\
            &=\text{HL}+\text{LH}.
        \end{aligned}
        \]
        For HL,
        \begin{equation}\label{5.24}
            \begin{aligned}
                \text{  HL} &\leq\sum_{l\in\mathbb{Z}\setminus\left\lbrace 0,k\right\rbrace, |k-l|\leq\frac{|k|}{2}} ||e^{\varepsilon_0\nu^{1/2}t}\hat{u}^1_l||_{L^2 L^\infty}||\hat{\theta}_{k-l}||_{L^\infty L^2}\\
                &\lesssim\sum_{l\in\mathbb{Z}\setminus \left\lbrace 0,k\right\rbrace ,|k-l|\leq\frac{|k|}{2}} \nu^{-\frac{1}{8}}|l|^{-\frac{3}{8}}E_lH_{k-l}\\
                &\lesssim\sum_{l\in\mathbb{Z}\setminus \left\lbrace 0,k\right\rbrace  ,|k-l|\leq\frac{|k|}{2}} \nu^{-\frac{1}{8}}|k|^{-\frac{3}{8}}E_lH_{k-l}\\
            \end{aligned}
        \end{equation}
        and
        \begin{equation}\label{5.25}
            \begin{aligned}
                \text{  LH} &\leq\sum_{l\in\mathbb{Z}\setminus \left\lbrace 0,k\right\rbrace  ,|k-l|\geq\frac{|k|}{2}} ||e^{\varepsilon_0\nu^{1/2}t}\hat{u}^1_l||_{L^\infty L^\infty}||\hat{\theta}_{k-l}||_{L^2 L^2}\\
                &\lesssim\sum_{l\in\mathbb{Z}\setminus\left\lbrace 0,k\right\rbrace ,|k-l|\geq\frac{|k|}{2}} (\nu|k-l|)^{-\frac{1}{4} \cdot \frac{1}{3}}\left(\nu^{\frac{1}{2}}|k-l|\right)^{-\frac{2}{3}}E_lH_{k-l}\\
                &\lesssim\sum_{l\in\mathbb{Z}\setminus \left\lbrace 0,k\right\rbrace ,|k-l|\geq\frac{|k|}{2}} \nu^{-\frac{5}{12}}|k|^{-\frac{3}{4}}E_lH_{k-l}.\\
            \end{aligned}
        \end{equation}
            And then, substituting \eqref{5.25} and \eqref{5.24} into \eqref{3.31}, we get
        \begin{equation}\label{5.26}
            \begin{aligned}
                ||e^{\varepsilon_0\nu^{1/2}t}g^1_k||_{L^2L^2}&\lesssim|k|^{-\frac{1}{2}}\nu^{-\frac{7}{24}}E_0H_k+\nu^{-\frac{1}{8}}|k|^{-\frac{3}{8}}E_kH_0\\
                &+\sum_{l\in\mathbb{Z}\setminus\left\lbrace 0,k\right\rbrace ,|k-l|\leq\frac{|k|}{2}} \nu^{-\frac{1}{8}}|k|^{-\frac{3}{8}}E_lH_{k-l}+\sum_{l\in\mathbb{Z}\setminus \left\lbrace 0,k\right\rbrace,|k-l|\geq\frac{|k|}{2}} \nu^{-\frac{5}{12}}|k|^{-\frac{3}{4}}E_lH_{k-l}.\\
            \end{aligned}
        \end{equation}
        The combination of equations  \eqref{eq:g2-estimate}, \eqref{eq:H-alpha-linear-from-theta} and \eqref{5.26} yields the result for $k\neq 0,$ and $\nu k^2\leq 1$
        \[
        \begin{aligned}
            &   H_k\leq||\hat{\theta}_{0}(k,y)||_{L^2}+C\left(  \nu^{-\frac{1}{4}}|k|^{\frac{3}{4}} ||e^{\varepsilon_0\nu^{1/2}t}g^1_k||_{L^2L^2}+\nu^{-\frac{1}{2}}||e^{\varepsilon_0\nu^{1/2}t}g^2_k||_{L^2L^2}\right)\\
            &\lesssim||\hat{\theta}_{0}(k,y)||_{L^2}+\nu^{-\frac{2}{3}}E_0H_k+\nu^{-\frac{1}{4}}\nu^{-\frac{1}{8}}|k|^{\frac{3}{8}}E_kH_0\\
            &+\sum_{l\in\mathbb{Z}\setminus\left\lbrace 0,k\right\rbrace ,|k-l|\leq\frac{|k|}{2}} \nu^{-\frac{1}{4}}\nu^{-\frac{1}{8}}|k|^{\frac{3}{8}}E_lH_{k-l}+\sum_{l\in\mathbb{Z}\setminus \left\lbrace 0,k\right\rbrace,|k-l|\geq\frac{|k|}{2}} \nu^{-\frac{2}{3}}E_lH_{k-l}+\nu^{-\frac{2}{3}}\sum_{l\in\mathbb{Z} } E_lH_{k-l}\\
            &\lesssim||\hat{\theta}_{0}(k,y)||_{L^2} +\nu^{-\frac{1}{8}}\nu^{-\frac{7}{16}}\sum_{l\in\mathbb{Z}\setminus\left\lbrace 0,k\right\rbrace ,|k-l|\leq\frac{|k|}{2}} E_lH_{k-l}\\
            &+\sum_{l\in\mathbb{Z}\setminus \left\lbrace 0,k\right\rbrace,|k-l|\geq\frac{|k|}{2}} \nu^{-\frac{2}{3}}E_lH_{k-l}+\nu^{-\frac{2}{3}}\sum_{l\in\mathbb{Z} } E_lH_{k-l}.\\
            &\lesssim|||\hat{\theta}_{0}(k,y)||_{L^2}+\sum_{l\in\mathbb{Z}} \nu^{-\frac{2}{3}}E_lH_{k-l}+\nu^{-\frac{1}{8}}\nu^{-\frac{7}{16}}\sum_{l\in\mathbb{Z}\setminus\left\lbrace 0,k\right\rbrace ,|k-l|\leq\frac{|k|}{2}} E_lH_{k-l}\\
        \end{aligned}
        \]

        Proof of \eqref{3.8}.
        For $\nu k^2\geq 1,$ form  Proposition \ref{prop:2.2}, we get
        \begin{equation}\label{5.27}
            \begin{aligned}
                &H_k=||e^{\varepsilon_0\nu^{1/2}t}\hat{\theta}_k||_{L^\infty L^2}+\nu^{\frac{1}{4}} k^{\frac{1}{4}}||e^{\varepsilon_0\nu^{1/2}t}\hat{\theta}_k||_{L^2L^2}+\nu^{\frac{1}{2}} |k|||e^{\varepsilon_0\nu^{1/2}t}\hat{\theta}_k||_{L^2L^2}\\
                &=\left( ||\hat{\theta}_k||_{L^\infty L^2}+|k|^{-\frac{3}{8}}(\nu k^2)^{\frac{1}{4}} ||e^{\varepsilon_0\nu^{1/2}t}\hat{\theta}_k||_{L^2L^2}+\nu^{\frac{1}{2}} |k|||e^{\varepsilon_0\nu^{1/2}t}\hat{\theta}_k||_{L^2L^2} \right)\\
                &\leq\left( ||\hat{\theta}_k||_{L^\infty L^2}+(\nu k^2)^{\frac{1}{2}} ||e^{\varepsilon_0\nu^{1/2}t}\hat{\theta}_k||_{L^2L^2}+\nu^{\frac{1}{2}} |k|||e^{\varepsilon_0\nu^{1/2}t}\hat{\theta}_k||_{L^2L^2}\right) \\
                &\lesssim\left( ||\hat{\theta}_0(k,y)||_{ L^2}+\nu ^{-\frac{1}{2}} \left( ||e^{\varepsilon_0\nu^{1/2}t}g^1_k||_{L^2L^2}+||e^{\varepsilon_0\nu^{1/2}t}g^2_k||_{L^2L^2}\right) \right). \\
            \end{aligned}
        \end{equation}
For the term
$$g^2_k(t,y)=\sum_{l\in\mathbb{Z}}\hat{u}^2_l\hat{\theta}_{k-l}(t,y),$$
noticing that $\bar{u}^2=0$, we also have that for $k\neq 0,$
        \begin{equation}\label{5.31}
            \begin{aligned}
                ||g^2_k||_{L^2L^2}\leq& \sum_{l\in\mathbb{Z},}\|e^{\varepsilon_0\nu^{1/2}t}\hat{u}^2_l\hat{\theta}_{k-l}\|_{L^2L^2}+\|e^{\varepsilon_0\nu^{1/2}t}\hat{u}^2_k\|_{L^2L^\infty}\|\bar{\theta}\|_{L^\infty L^2}\\
                \leq
                &\sum_{l\in\mathbb{Z} ,l\neq \left\lbrace 0,k\right\rbrace } ||\hat{u}^2_l||_{L^\infty L^\infty}||e^{\varepsilon_0\nu^{1/2}t}\hat{\theta}_{k-l}||_{L^2 L^2}+ ||\hat{u}^2_k||_{L^2 L^2}^{\frac{1}{2}}||e^{\varepsilon_0\nu^{1/2}t}k\hat{u}^1_k\|_{L^2L^2}^{\frac{1}{2}}H_0\\
                \leq&\sum_{l\in\mathbb{Z} } ||\hat{u}^2_l||_{L^\infty L^2}^{\frac{1}{2}}||\partial_y\hat{u}^2_l\|_{L^\infty L^2}^{\frac{1}{2}}||e^{\varepsilon_0\nu^{1/2}t}\hat{\theta}_{k-l}||_{L^2 L^2}+|k|^{-\frac{1}{2}}E_kH_0\\
                \leq&\sum_{l\in\mathbb{Z} } ||\hat{u}^2_l||_{L^\infty L^2}^{\frac{1}{2}}||l\hat{u}^1_l\|_{L^\infty L^2}^{\frac{1}{2}}||e^{\varepsilon_0\nu^{1/2}t}\hat{\theta}_{k-l}||_{L^2 L^2}+|k|^{-\frac{1}{2}}E_kH_0\\
                \leq&\sum_{l\in\mathbb{Z} }|l|^{-\frac{1}{2}}\nu^{-\frac{1}{4}} |k-l|^{-\frac{1}{4}}E_lH_{k-l}+|k|^{-\frac{1}{2}}E_kH_0.\\
                \leq&\sum_{l\in\mathbb{Z} }\nu^{-\frac{1}{4}} |k|^{-\frac{1}{4}}E_lH_{k-l}+|k|^{-\frac{1}{2}}E_kH_0.\\
            \end{aligned}
        \end{equation}
        Combining \eqref{5.27}, \eqref{5.31},  for $\nu k^2\geq 1$, we get
 \[
 \begin{aligned}
 H_k&
 \leq\left( ||\hat{\theta}_0(k,y)||_{ L^2}+\nu ^{-\frac{1}{2}} \left( ||e^{\varepsilon_0\nu^{1/2}t}g^1_k||_{L^2L^2}+||e^{\varepsilon_0\nu^{1/2}t}g^2_k||_{L^2L^2}\right) \right), \\
 &\lesssim||\hat{\theta}_{0}(k,y)||_{L^2}+|k|^{-\frac{3}{4}}\nu^{-\frac{11}{12}}E_0H_k+\nu^{-\frac{1}{2}}\nu^{-\frac{1}{8}}|k|^{-\frac{3}{8}} E_kH_0+\nu^{-\frac{1}{2}}|k|^{-\frac{1}{2}} E_kH_0\\
 &+\sum_{l\in\mathbb{Z}\setminus \left\lbrace 0,k\right\rbrace} \nu^{-\frac{3}{4}}|k|^{-\frac{1}{4}}E_lH_{k-l}+|k|^{-\frac{3}{4}}\nu^{-\frac{11}{12}}\sum_{l\in\mathbb{Z} } E_lH_{k-l}\\
 &\lesssim||\hat{\theta}_{0}(k,y)||_{L^2}+\nu^{-\frac{5}{8}}\sum_{l\in\mathbb{Z}\setminus\{0\}} E_lH_{k-l}+C\nu^{-\frac{5}{8}}\sum_{l\in\mathbb{Z}}E_lH_{k-l}+C\nu^{-\frac{13}{24}}\sum_{l\in\mathbb{Z}}E_lH_{k-l}\\
 &\lesssim||\hat{\theta}_{0}(k,y)||_{L^2}+C\nu^{-\frac{5}{8}}\sum_{l\in\mathbb{Z}}E_lH_{k-l}+C\nu^{-\frac{13}{24}}\sum_{l\in\mathbb{Z}}E_lH_{k-l}\\
 \end{aligned}
 \]
 This completes the proof of Proposition \ref{pro3.1}.
 \end{proof}
 Theorem \ref{thm1.2} is now proven. By utilizing \eqref{3.4} and \eqref{3.5}, we derive
                    \begin{equation}\label{eq:3.40}
\sum_{k\in\mathbb Z}E_k
\lesssim
\sum_{k\in\mathbb Z}\|\omega_k^{\rm in}\|_{H_k^4}
+
\nu^{-2/3}
\sum_{k\in\mathbb Z}
\sum_{l\in\mathbb Z}
 E_l E_{k-l}+\sum_{k\in\mathbb Z}\nu^{-1/2}\nu^{-1/4}H_k.
\end{equation}
Recalling the fact
$$\sum_{k\in\mathbb{Z}}H_k=H_0+\sum_{k\in\mathbb{Z}\setminus\{0\}, \nu k^2\leq 1}H_k+\sum_{k\in\mathbb{Z}\setminus\{0\},\nu k^2>1}H_k,$$
combining \eqref{3.6}, \eqref{3.7} and \eqref{3.8}, we can deduce
\begin{equation}\label{5.33}
\begin{aligned}
&\sum_{k\in\mathbb{Z}}H_k\lesssim||\bar{\theta}_{0}||_{L^2}
+\sum_{k\in\mathbb{Z}\setminus\{0\} }||\hat{\theta}_{0}(k,y)||_{L^2}
+\nu^{-\frac{5}{8}}\sum_{k\in\mathbb{Z} }\sum_{l\in\mathbb{Z} } E_lH_{k-l}\\
&+\nu^{-\frac{1}{8}}\nu^{-\frac{7}{16}}\sum_{k\in\mathbb{Z}\setminus\{0\},
\nu k^2\leq 1}\sum_{l\in\mathbb{Z}} E_lH_{k-l}
+\nu^{-\frac{13}{24}}\sum_{k\in\mathbb{Z}\setminus\{0\},\nu k^2>1}\sum_{l\in\mathbb{Z}}E_lH_{k-l}.\\
\end{aligned}
\end{equation}
On the other hand, if
$$||u_{0}||_{H^{5+}}\leq c_0\min\left\lbrace
\nu,\nu\right\rbrace ^{\frac{2}{3}},\quad ||\theta_{0}||_{H^1}\leq
c_1\min\left\lbrace \nu,\nu\right\rbrace ^{{\frac{17}{12}}},$$ we
have the ability to obtain
$$\sum_{k\in\mathbb Z}\|\omega_k^{\rm in}\|_{H_k^4}\leq Cc_0\min\left\lbrace \nu,\nu\right\rbrace ^{\frac{2}{3}},$$
and
$$||\bar{\theta}_{0}||_{L^2}+\sum_{k\in\mathbb{Z}\setminus\{0\}}||\hat{\theta}_{0}(k,y)||_{L^2}
\leq Cc_1\min\left\lbrace \nu,\nu\right\rbrace ^{\frac{17}{12}}.$$
Thus,  selecting $c_0,c_1$ suitably small, by bootstrap arguments,
and from \eqref{eq:3.40} and \eqref{5.33}, we can infer that
$$\sum_{k\in\mathbb{Z}}E_k\leq Cc_0\min\left\lbrace \nu,\nu\right\rbrace ^{\frac{2}{3}},\quad\sum_{k\in\mathbb{Z}}H_k
\leq Cc_1\min\left\lbrace \nu,\nu\right\rbrace ^{\frac{17}{12}}.$$
This completes the proof of Theorem \ref{thm1.2}. \

 {\bf Acknowledgements:}  The research is supported in part
by the National Natural Science Foundation of China (Grant Nos.
12331007,12526581 and 12671260).

\textbf{Conflict of Interest}\quad The authors declared that they
have no conflict of interest.

\end{document}